\documentclass[twoside,11pt]{article}
\usepackage{a4wide,latexsym,amsfonts,amssymb,exscale,epsfig}
\usepackage[centertags,sumlimits,intlimits,namelimits,reqno]{amsmath}

\usepackage{amsthm}
\theoremstyle{definition}
\newtheorem{theorem}{Theorem}[section]

\newtheorem{corollary}[theorem]{Corollary}
\newtheorem{proposition}[theorem]{Proposition}
\newtheorem{definition}[theorem]{Definition}
\newtheorem{remark}[theorem]{Remark}
\newtheorem{example}[theorem]{Example}

\usepackage{hyperref}
\usepackage{cite}

\usepackage{bbm}
\def\C{{\mathbbm C}}
\def\N{{\mathbbm N}}

\def\ie{{\sl i.e.\/}}

\def\etc{{\sl etc.\/}}
\def\cf{{\sl c.f.\/}}
\let\phi=\varphi
\let\theta=\vartheta
\let\epsilon=\varepsilon

\def\dim{\mathop{\rm dim}\nolimits}

\def\id{\mathop{\rm id}\nolimits}
\def\Span{\mathop{\rm span}\nolimits}

\def\GL{\mathop{\rm GL}\nolimits}
\def\Hom{\mathop{\rm Hom}\nolimits}
\def\Aut{\mathop{\rm Aut}\nolimits}

\def\ker{\mathop{\rm ker}\nolimits}

\def\lim{\mathop{\rm lim}\limits}

\let\hat=\widehat
\let\tilde=\widetilde

\newcommand{\color}[2][c]{}

\numberwithin{equation}{section}

\makeatletter

\newfont{\@aidxte}{cmsy10}
\newfont{\@aidxel}{cmsy10 scaled 1095}
\newfont{\@aidxtw}{cmsy10 scaled 1200}
\newlength\@aidxtexvi
\newlength\@aidxtexvii
\newlength\@aidxelxvi
\newlength\@aidxelxvii
\newlength\@aidxtwxvi
\newlength\@aidxtwxvii
\newcommand{\alignidx}[1]{%
  \@aidxtexvi=\fontdimen16\@aidxte
  \@aidxtexvii=\fontdimen17\@aidxte
  \@aidxelxvi=\fontdimen16\@aidxel
  \@aidxelxvii=\fontdimen17\@aidxel
  \@aidxtwxvi=\fontdimen16\@aidxtw
  \@aidxtwxvii=\fontdimen17\@aidxtw
    {\mbox{$%
    \fontdimen16\@aidxte=2.9pt
    \fontdimen17\@aidxte=2.9pt
    \fontdimen16\@aidxel=3.1pt
    \fontdimen17\@aidxel=3.1pt
    \fontdimen16\@aidxtw=3.3pt
    \fontdimen17\@aidxtw=3.3pt
    #1$}}%
    \fontdimen16\@aidxte=\@aidxtexvi
    \fontdimen17\@aidxte=\@aidxtexvii
    \fontdimen16\@aidxel=\@aidxelxvi
    \fontdimen17\@aidxel=\@aidxelxvii
    \fontdimen16\@aidxtw=\@aidxtwxvi
    \fontdimen17\@aidxtw=\@aidxtwxvii}

\makeatother

\newenvironment{myenumerate}{%
  \begin{enumerate}
  \setlength{\partopsep}{0pt}
  \setlength{\parskip}{0pt}}{\end{enumerate}}

\newenvironment{myitemize}{%
  \begin{itemize}
  \setlength{\itemsep}{0pt}
  \setlength{\parskip}{0pt}}{\end{itemize}}

\def\nn{\notag}
\def\del{\partial}
\def\emph#1{{\sl #1\/}}

\def\K{k}

\def\op{\mathrm{op}}

\def\Set{\mathbf{Set}}
\def\compHaus{\mathbf{compHaus}}
\def\comUnCstAlg{\mathbf{comUnC^\ast\!Alg}}
\def\Grp{\mathbf{Grp}}
\def\CoGrp{\mathbf{CoGrp}}
\def\fGrp{\mathbf{fGrp}}
\def\Cat{\mathbf{Cat}}
\def\CoCat{\mathbf{CoCat}}
\def\fcCat{\mathbf{Lex}}
\def\fpCat{\mathbf{FP}}
\def\TCat{\mathbf{2Cat}}
\def\TGrp{\mathbf{2Grp}}
\def\XMod{\mathbf{XMod}}
\def\fTGrp{\mathbf{f2Grp}}
\def\compTopGrp{\mathbf{compTopGrp}}
\def\compTopTGrp{\mathbf{compTop2Grp}}
\def\Vect{\mathbf{Vect}}
\def\fdVect{\mathbf{fdVect}}
\def\Alg{\mathbf{Alg}}
\def\Co{\mathbf{CoAlg}}
\def\Bi{\mathbf{BiAlg}}
\def\Hopf{\mathbf{HopfAlg}}
\def\cocHopf{\mathbf{cocHopfAlg}}
\def\cocBi{\mathbf{cocBiAlg}}
\def\cocCo{\mathbf{cocCoAlg}}
\def\cocTriAlg{\mathbf{cocTriAlg}}
\def\Tri{\mathbf{TriAlg}}
\def\CoTri{\mathbf{CoTriAlg}}
\def\comCoTri{\mathbf{comCoTriAlg}}
\def\comHopf{\mathbf{comHopfAlg}}
\def\comHopfCstAlg{\mathbf{comHopfC^\ast\!Alg}}
\def\comCstCoTri{\mathbf{comC^\ast\!CoTriAlg}}
\def\comBi{\mathbf{comBiAlg}}
\def\comAlg{\mathbf{comAlg}}
\def\Th{\mathbf{Th}}
\def\Mod{\mathbf{Mod}}
\def\coend{\mathbf{coend}}
\def\comod{\mathbf{comod}}
\def\symHopf{\mathbf{symHopfCat}}
\def\HopfCat{\mathbf{HopfCat}}
\def\Rep{\mathbf{Rep}}

\def\pullback#1#2#3{%
  \,\mbox{\raisebox{-.8ex}{$\scriptstyle #1$}}%
  \!\prod\!
  \mbox{\raisebox{-.8ex}{$\scriptstyle #3$}}\,}%
\def\dpullback#1#2#3{\prod\limits_{#1\,\,#3}}
\def\pushout#1#2#3{%
  \,\mbox{\raisebox{-.8ex}{$\scriptstyle #1$}}%
  \!\coprod\!
  \mbox{\raisebox{-.8ex}{$\scriptstyle #3$}}\,}%

\def\Uepsilon{\underline{\epsilon}}
\def\UDelta{\underline{\Delta}}
\def\US{\underline{S}}

\def\eq{\mathrm{eq}}

\def\acknowledgments{\section*{Acknowledgments}}%
\usepackage[ps,matrix,arrow,curve]{xy}

\newcommand{\hpeprint}[1]{%
  \href{http://arXiv.org/abs/#1}{\texttt{#1}}}%
\newcommand{\hpmathsci}[1]{%
  \href{http://www.ams.org/mathscinet-getitem?mr=#1}{\texttt{MR #1}}}%

\begin{document}
%

\title{
\hbox to\hsize{\hfill\small DAMTP-2004-135}\vspace{-4mm}
\hbox to\hsize{\hfill\small ESI 1603 (2005)}
$2$-Groups, trialgebras and their Hopf categories\\
       of representations}
\author{Hendryk Pfeiffer\thanks{E-mail: \texttt{H.Pfeiffer@damtp.cam.ac.uk}}}
\date{\small{DAMTP, Wilberforce Road, Cambridge CB3 0WA, United Kingdom}\\
  \small{Emmanuel College, St Andrew's Street, Cambridge CB2 3AP, United
  Kingdom}\\[1ex]
  September 27, 2006
}

\maketitle

%
\begin{abstract}
%

A strict $2$-group is a $2$-category with one object in which all
morphisms and all $2$-morphisms have inverses. $2$-Groups have been
studied in the context of homotopy theory, higher gauge theory and
Topological Quantum Field Theory (TQFT). In the present article, we
develop the notions of trialgebra and cotrialgebra, generalizations of
Hopf algebras with two multiplications and one comultiplication or
vice versa, and the notion of Hopf categories, generalizations of
monoidal categories with an additional functorial comultiplication. We
show that each strict $2$-group has a `group algebra' which is a
cocommutative trialgebra, and that each strict finite $2$-group has a
`function algebra' which is a commutative cotrialgebra. Each such
commutative cotrialgebra gives rise to a symmetric Hopf category of
corepresentations. In the semisimple case, this Hopf category is a
$2$-vector space according to Kapranov and Voevodsky. We also show
that strict compact topological $2$-groups are characterized by their
$C^\ast$-cotrialgebras of `complex-valued functions', generalizing the
Gel'fand representation, and that commutative cotrialgebras are
characterized by their symmetric Hopf categories of corepresentations,
generalizing Tannaka--Kre\v\i n reconstruction. Technically, all these
results are obtained using ideas from functorial semantics, by
studying models of the essentially algebraic theory of categories in
various base categories of familiar algebraic structures and the
functors that describe the relationships between them.
\end{abstract}

%
\section{Introduction}
%

A strict $2$-group is an internal category in the category of
groups. Strict $2$-groups can also be characterized as $2$-categories
with one object in which all morphisms and all $2$-morphisms have
inverses, \ie\ as (strict) $2$-groupoids with one object. The notion
of a strict $2$-group can therefore be viewed as a higher-dimensional
generalization of the notion of a group because the set of
$2$-morphisms of a $2$-group has got \emph{two} multiplication
operations, a horizontal one and a vertical one\footnote{The prefix
`$2$-' obviously refers to the higher-dimensional nature rather than
to the order of the group.}. Strict $2$-groups can be constructed from
Whitehead's crossed modules, and so there exist plenty of examples.

Starting from the theory of groups, one can develop the notion of
cocommutative Hopf algebras which arise as group algebras, the notion
of commutative Hopf algebras which appear as algebras of functions on
groups, and the notion of symmetric monoidal categories which arise as
the representation categories of groups. Compact topological groups
are characterized by their commutative Hopf $C^\ast$-algebras of
continuous complex-valued functions (Gel'fand
representation). Commutative Hopf algebras are characterized by their
rigid symmetric monoidal categories of finite-dimensional comodules
(Tannaka--Kre\v\i n reconstruction). This `commutative' theory forms
the basic framework that is required before one can develop the theory
of quantum groups.

What is the analogue of the preceding paragraph if one systematically
replaces the word `group' by `strict $2$-group'? In particular, what
is a good definition of `group algebra', `function algebra' and
`representation category' for a strict $2$-group? The sought-after
definitions can be considered successful if they allow us to retain
analogues of the most important theorems that are familiar from
groups. It is the purpose of the present article to present concise
definitions of the relevant structures and to establish
generalizations from groups to strict $2$-groups of the theorems
mentioned above, namely on Gel'fand representation and on
Tannaka--Kre\v\i n reconstruction. We hope that our definitions and
results will prove useful in order to investigate whether there exist
structures that can play the role of `quantum $2$-groups', but this
question lies beyond the scope of the present article.

\subsection{Higher-dimensional algebra}

With the step from groups to strict $2$-groups, we enter the realm of
higher-dimensional algebra. Higher-dimensional algebraic structures
have appeared in various areas of mathematics and mathematical
physics. A prime example is the higher-dimensional group theory
programme of Brown~\cite{Br82}, generalizing groups and groupoids to
double groupoids and further on, in order to obtain a hierarchy of
algebraic structures. The construction of these algebraic structures
is inspired by problems in homotopy theory where algebraic structures
at a some level of the hierarchy are related to topological features
that appear in the corresponding dimension.

In order to construct topological quantum field theories (TQFTs),
Crane has introduced the concept of \emph{categorification}, see, for
example~\cite{CrFr94,BaDo98}. Categorification can be viewed as a
systematic replacement of familiar algebraic structures that are
modelled on sets by analogues that are rather modelled on categories,
$2$-categories, and so on. Categorification often serves as a guiding
principle in order to find suitable definitions of algebraic
structures at some higher level starting from the known definitions at
a lower level.

Some examples of higher-dimensional algebraic structures that are
relevant in the context of the present article, are the following.
\begin{myitemize}
\item
  Three-dimensional TQFTs can be constructed from Hopf
  algebras~\cite{Ku91,ChFu94}. In order to find four-dimensional
  TQFTs, Crane and Frenkel~\cite{CrFr94} have introduced the notion of
  a \emph{Hopf category}, a categorification of the notion of a Hopf
  algebra. Roughly speaking, this is a monoidal category with an
  additional functorial comultiplication.
\item
  Crane and Frenkel~\cite{CrFr94} have also speculated about
  \emph{trialgebras}, vector spaces with three mutually compatible
  linear operations: two multiplications and one comultiplication or
  vice versa. Hopf categories are thought to appear as the
  representation categories of these trialgebras in analogy to
  monoidal categories which appear as the representation categories of
  Hopf algebras.
\item
  Kapranov and Voevodsky~\cite{KaVo91,KaVo94} have introduced braided
  monoidal $2$-ca\-te\-go\-ries and $2$-vector spaces, a categorified
  notion of vector spaces, and shown that they are related to the
  Zamolodchikov tetrahedron equation. In the context of integrable
  systems in mathematical physics, this equation is thought to be the
  generalization of the integrability condition to three
  dimensions. In two dimensions, the integrability condition is the
  famous Yang--Baxter equation.
\item
  Grosse and Schlesinger have constructed examples of
  trialgebras~\cite{GrSc00a,GrSc00b} in the spirit of Crane--Frenkel
  and explained how they are related to integrability in $2+1$
  dimensions~\cite{GrSc01}.
\item
  Baez, Lauda, Crans, Bartels, Schreiber and the
  author~\cite{Ba02,BaLa04,BaCr04,Pf03,Ba04,BaSc04} have used
  $2$-groups in order to find generalizations of fibre bundles and of
  gauge theory. Yetter~\cite{Ye93} has used $2$-groups in order to
  construct novel TQFTs, generalizing the TQFTs that are constructed
  from the gauge theories of flat connections on a principal
  $G$-bundle where $G$ is an (ordinary) group. One can verify
  that~\cite{Ye93} is a special application of the generalized gauge
  theory of~\cite{Pf03}.
\end{myitemize}

All these constructions have a common underlying theme: the procedure
of categorification on the algebraic side and an increase in dimension
on the topological side. Although it is plausible to conjecture that
all these higher-dimensional algebraic structures are related, the
developments mentioned above have so far been largely independent. The
present article with its programme of finding the `group algebras',
`function algebras' and `representation categories' of strict
$2$-groups provides relationships between these structures. In the
remainder of the introduction, we sketch our main results and explain
how these can be seen in the context of the existing literature.

\subsection{$2$-Groups}

A strict $2$-group $(G_0,G_1,s,t,\imath,\circ)$ is an internal
category in the category of groups. It therefore consists of groups
$G_0$ (\emph{group of objects}) and $G_1$ (\emph{group of morphisms})
and of group homomorphisms $s\colon G_1\to G_0$ (\emph{source}),
$t\colon G_1\to G_0$ (\emph{target}), $\imath\colon G_0\to G_1$
(\emph{identity}) and $\circ\colon G_1\times_{G_0} G_1\to G_1$
(\emph{composition}) where $G_1\times_{G_0}G_1=\{\,(g,g^\prime)\in
G_1\times G_1\colon\, t(g)=s(g^\prime)\,\}$ denotes the group of
composable morphisms, subject to conditions which are given in detail
in Section~\ref{sect_internalization}.

\subsection{$2$-Groups, trialgebras and cotrialgebras}
\label{sect_introtrialg}

Passing from a group $G$ to its group algebra $\K[G]$ is described by
a functor $\K[-]$ from the category of groups to the category of
cocommutative Hopf algebras over some field $\K$. Applying this
functor systematically to all groups and group homomorphisms in the
definition of a strict $2$-group yields a structure
$(H_0,H_1,s^\prime,t^\prime,\imath^\prime,\circ^\prime)=(\K[G_0],\K[G_1],\K[s],\K[t],\K[\imath],\K[\circ])$
which consists of cocommutative Hopf algebras $H_0$ and $H_1$ and
bialgebra homomorphisms $s^\prime$, $t^\prime$, $\imath^\prime$ and
$\circ^\prime$. What sort of algebraic structure is this?

Since the functor $\K[-]$ preserves all finite limits,
$(H_0,H_1,s^\prime,t^\prime,\imath^\prime,\circ^\prime)$ is an
internal category in the category of cocommutative Hopf algebras over
$\K$ (Section~\ref{sect_cochopflimit}). Unfolding this definition
(Section~\ref{sect_coctrialgdef}), we see that $H_1=\K[G_1]$ is a
$\K$-vector space with three mutually compatible linear operations:
two multiplications and one comultiplication. We thus call
$(H_0,H_1,s^\prime,t^\prime,\imath^\prime,\circ^\prime)$ a
\emph{cocommutative trialgebra}. This yields a concise definition for
the structure conjectured by Crane and Frenkel~\cite{CrFr94} as an
internal category in the category of cocommutative Hopf algebras and
allows us to find a large class of examples.

Similarly to $\K[-]$, there is a functor $\K(-)$ from the category of
finite groups to the category of commutative Hopf algebras which sends
every finite group $G$ to the algebra of $\K$-valued functions on
$G$. By a procedure similar to that for $\K[-]$, we arrive at the
definition of a \emph{commutative cotrialgebra} and a method for
constructing examples from strict finite $2$-groups
(Section~\ref{sect_cotrialg}). We generalize Gel'fand representation
theory to $2$-groups and prove that any strict compact topological
$2$-group gives rise to a commutative $C^\ast$-cotrialgebra and,
conversely, that the compact topological $2$-group can be
reconstructed from this $C^\ast$-cotrialgebra
(Section~\ref{sect_top2grp}).

\subsection{Cotrialgebras and Hopf categories}

For the relationship with Hopf categories and $2$-vector spaces, we
refer to the scenario of Crane and Frenkel~\cite{CrFr94} in greater
detail.

When one tries to construct topological invariants or TQFTs for
combinatorial manifolds, one defines so-called state
sums~\cite{FuHo94,Ku91,ChFu94,TuVi92,BaWe96}. It turns out that for
every dimension of the manifolds, there exist preferred algebraic
structures that guarantee the consistency and triangulation
independence of these state sums.

Whereas $(1+1)$-dimensional TQFTs can be constructed from suitable
associative algebras~\cite{FuHo94}, there are \emph{two} alternative
algebraic structures in order to construct $(2+1)$-dimensional TQFTs:
suitable Hopf algebras~\cite{ChFu94,Ku91} or suitable monoidal
categories~\cite{TuVi92,BaWe96}. Both types of structure are
related~\cite{BaWe95}: the category of corepresentations of a Hopf
algebra is a monoidal category and, conversely, under some conditions,
the Hopf algebra can be Tannaka--Kre\v\i n reconstructed from this
category.

According to Crane and Frenkel~\cite{CrFr94}, there ought to be
\emph{three} alternative algebraic structures in order to construct
$(3+1)$-dimensional TQFTs. These have been suggestively termed
\emph{trialgebras}, \emph{Hopf categories} and \emph{monoidal
$2$-categories}. All three types of structures are thought to be
related: representing one product of a trialgebra should yield a Hopf
category as the category of representations of the trialgebra, and
then representing the second product should give a monoidal
$2$-category as the $2$-category of representations of the Hopf
category. Combining both steps, \ie\ representing both products at
once, should give a monoidal $2$-category as the category of
representations of the trialgebra.

Some aspects of the Crane--Frenkel scenario have already been analyzed
in greater detail.
\begin{myitemize}
\item
  Mackaay~\cite{Ma99} has given a precise definition of suitable
  monoidal $2$-categories and has shown that one can define an
  invariant of combinatorial $4$-manifolds from it. So far, only few
  examples of these monoidal $2$-categories have been constructed all
  of which are thought to give homotopy invariants.
\item
  Neuchl~\cite{Ne99} has studied Hopf categories and their
  representations on certain monoidal $2$-categories. So far, it is
  open whether one can find a good class of Hopf categories and the
  corresponding monoidal $2$-categories in such a way that both
  structures are related by a generalization of Tannaka--Kre\v\i n
  duality.
\item
  Crane and Frenkel~\cite{CrFr94} have presented examples of Hopf
  categories and proposed a $4$-manifold invariant based on Hopf
  categories. Carter, Kauffman and Saito~\cite{CaKa99} have studied
  this invariant for some Hopf categories. Again, with the rather
  limited set of examples of Hopf categories, all invariants studied
  so far, are homotopy invariants.
\end{myitemize}

In the present article, we extend the study of the Crane--Frenkel
scenario and show in addition that the corepresentations of a
commutative cotrialgebra form a symmetric Hopf category and,
conversely, we prove a generalization of Tannaka--Kre\v\i n duality in
order to recover the cotrialgebra from its Hopf category of
corepresentations.

Combining this with the result sketched in
Section~\ref{sect_introtrialg} that each strict compact topological
$2$-group gives rise to a commutative cotrialgebra, this shows that
each compact topological $2$-group has got a symmetric Hopf category
of finite-dimensional continuous unitary representations and that the
$2$-group can be Tannaka--Kre\v\i n reconstructed from this Hopf
category.

In order to make these theorems possible, our definition of a
symmetric Hopf category deviates from those used in the literature by
Crane--Frenkel and by Neuchl~\cite{CrFr94,Ne99} in a subtle way
(Section~\ref{sect_hopfcat}), in particular: (1) the functorial
comultiplication maps into a pushout rather than into the external
tensor product; (2) the comultiplication is already uniquely specified
as soon as certain other data are given; (3) it automatically
possesses a functorial antipode. Our definition is designed in such a
way that Tannaka-Kre\v\i n reconstruction generalizes to $2$-groups
and that there is the notion of the representation category of a
strict $2$-group which forms a symmetric Hopf category. In future work
on TQFTs, it may be useful to employ a generalization of our
definition of a Hopf category as opposed to one in which the
comultiplication maps into the external tensor product. Imagine one
tried to invent the concept of a monoidal category, but one was
unlucky and proposed a definition that did not include the appropriate
representation categories of groups as examples.

Our Hopf categories of corepresentations of cotrialgebras form a
special case of $2$-vector spaces according to Kapranov and
Voevodsky~\cite{KaVo91} as soon as the cotrialgebra is
cosemisimple. This includes in particular the cotrialgebras of
complex-valued representative functions of strict compact topological
$2$-groups and provides us with a generalization of the Peter--Weyl
decomposition from compact topological groups to strict compact
topological $2$-groups.

Other authors have explored alternative strategies for representing
$2$-groups on $2$-vector spaces. Barrett and Mackaay~\cite{BaMa04} and
Elgueta~\cite{El04}, for example, employ the $2$-vector spaces of
Kapranov and Voevodsky~\cite{KaVo91} (semisimple Abelian categories),
Crane and Yetter~\cite{CrYe03b,Ye05} use a measure theoretic
refinement of these in order to include more interesting examples,
whereas Forrester-Barker~\cite{Fo03} uses the $2$-vector spaces of
Baez and Crans~\cite{BaCr04} (internal categories in the category of
vector spaces or, equivalently, $2$-term chain complexes of vector
spaces). All these authors exploit the fact that one can associate
with each $2$-group a $2$-category with one object, just as one can
associate with each group a category with one object. A representation
of the $2$-group is then defined to be a $2$-functor from this
$2$-category to the $2$-category of $2$-vector spaces, generalizing
the fact that a representation of an ordinary group is a functor from
the corresponding category with one object to the category of vector
spaces. Whereas the key construction of these approaches is a
$2$-functor from \emph{a single} $2$-group to the $2$-category of
$2$-vector spaces, in the present article we employ a $2$-functor from
the $2$-category of \emph{all} $2$-groups to the $2$-category of
$2$-vector spaces (in fact Hopf categories). In our case, the Hopf
category therefore plays the role of the \emph{representation
category} of the $2$-group rather than that of an \emph{individual
representation}. Our result on Tannaka--Kre\v\i n reconstruction
confirms that this is a useful approach, too.

\subsection{Category theoretic techniques}
\label{sect_categories}

In Section~\ref{sect_introtrialg}, we have sketched how to apply the
finite-limit preserving functors $\K[-]$ and $\K(-)$ to all objects
and morphisms in the definition of a $2$-group in order to obtain the
definitions of cocommutative trialgebras and commutative
cotrialgebras. The application of such finite-limit preserving
functors is the main theme of the present article.

In the more sophisticated language of functorial
semantics~\cite{La63}, the key technique exploited in the present
article is the study of models of the essentially algebraic theory of
categories in various interesting base categories of familiar
algebraic structures such as groups, compact topological groups,
commutative or cocommutative Hopf algebras or symmetric monoidal
categories and to develop their relationships.

\subsection{Limitations}

It is beyond the scope of the present work to find the `most generic'
examples of these novel structures or to attempt a
classification. Since we use only an abstract machinery in order to
construct higher-dimensional algebraic structures from conventional
ones and since it is natural to expect that the higher-dimensional
structures are substantially richer than the conventional ones, our
tools are not sufficient in order to survey the entire new
territory. What we can achieve is to find a first path through the new
structures that includes the strict and the commutative or symmetric
special cases. We expect two sorts of generalizations beyond the
present work. Firstly, the use of weak $2$-groups~\cite{BaLa04} rather
than strict ones and a similar weakening of the notions of trialgebra
and Hopf category. Weak $2$-groups are modelled on bicategories
whereas strict ones are modelled on $2$-categories. Secondly,
classical constructions starting from ordinary groups alone yield only
commutative or cocommutative Hopf algebras, but not actual quantum
groups. Quantum groups rather involve a novel idea beyond the
classical machinery such as Drinfeld's quantum double construction. At
the higher-dimensional level, we expect (at least) the same
limitations. We can construct infinite families of examples of
cocommutative trialgebras, of commutative cotrialgebras and of
symmetric Hopf categories in a systematic way, but our method does not
deliver any actual `quantum $2$-groups' yet.

Ultimately, one will need applications of the novel structures in
topology and in mathematical physics in order to confirm whether
sufficiently generic examples of `quantum $2$-groups' have been
found. Test cases are firstly the systematic construction of
non-trivial integrable systems in three dimensions (as opposed to the
usual two dimensions); secondly the construction of four-dimensional
TQFTs which give rise to invariants of differentiable four-manifolds
that are finer than just homotopy invariants; thirdly to overcome
problems in the spin foam approach to the quantization of general
relativity, see, for example~\cite{CrYe03a}.

\subsection{Outline}

This article is structured as follows. Section~\ref{sect_prelim} fixes
the notation. We summarize the definition of internal categories in
finitely complete base categories and the definition of strict
$2$-groups and some important examples of these. We define and
construct examples of cocommutative trialgebras in
Section~\ref{sect_coctrialg}. Commutative cotrialgebras are defined
and constructed in Section~\ref{sect_comcotrialg}. We introduce
symmetric Hopf categories in Section~\ref{sect_hopfcat} and develop
the corepresentation theory of cotrialgebras and the representation
theory of compact topological $2$-groups. More generic trialgebras and
Hopf categories are briefly mentioned in
Section~\ref{sect_gentrialg}. 

%
\section{Preliminaries}
%
\label{sect_prelim}

In this section, we fix some notation, we briefly review the technique
of internalization in order to construct the $2$-categories of
internal categories in familiar base categories, and we recall the
definition of strict $2$-groups as internal categories in the category
of groups.

\subsection{Notation}
\label{sect_notation}

We fix a field $\K$. Our notation for some standard categories is as
follows: $\Set$ (sets and maps), $\compHaus$ (compact Hausdorff spaces
and continuous maps), $\comUnCstAlg$ (commutative unital
$C^\ast$-algebras with unital $\ast$-homomorphisms), $\Grp$ (groups
and group homomorphisms), $\compTopGrp$ (compact topological groups),
$\Vect_\K$ ($\K$-vector spaces with $\K$-linear maps), $\fdVect_\K$
(finite-dimensional $\K$-vector spaces with $\K$-linear maps),
$\Alg_\K$ (associative unital algebras over $\K$ and their
homomorphisms), $\Co_\K$ (coalgebras over $\K$ and their
homomorphisms), $\Bi_\K$ (bialgebras) and $\Hopf_\K$ (Hopf
algebras). We use the prefix $\mathbf{f}$ for `finite', $\mathbf{fd}$
for `finite-dimensional', $\mathbf{com}$ for `commutative' and
$\mathbf{coc}$ for `cocommutative'.

Let $\mathcal{C}$ be some category. For each object $A$ of
$\mathcal{C}$, we denote its identity morphism by $\id_A\colon A\to
A$. Composition of morphisms $f\colon A\to B$ and $g\colon B\to C$ in
$\mathcal{C}$ is denoted by juxtaposition $fg\colon A\to C$ and is
read from left to right.

For objects $A$ and $B$ of $\mathcal{C}$, we denote their product (if
it exists) by $A\prod B$ and by $p_1\colon A\prod B\to A$ and
$p_2\colon A\prod B\to B$ the associated morphisms of its limiting
cone. Similarly, the coproduct (if it exists) is called $A\coprod B$,
and $\imath_1\colon A\to A\coprod B$ and $\imath_2\colon B\to A\coprod
B$ the morphisms of its colimiting cone.

Let $s\colon A\to C$ and $t\colon B\to C$ be morphisms of
$\mathcal{C}$. We denote by 
\begin{equation}
  A\textstyle\pullback{s}{C}{t}B:=\lim_{\leftarrow}(\xymatrix@1{A\ar[r]^s&C&B\ar[l]_t}),
\end{equation}
the pullback object (if it exists) and by $p_1\colon
A\pullback{s}{C}{t}B\to A$ and $p_2\colon A\pullback{s}{C}{t}B\to B$
its limiting cone, \ie,
\begin{equation}
\begin{aligned}  
\xymatrix{
  A\pullback{s}{C}{t}B\ar[rr]^{p_2}\ar[dd]_{p_1}&&B\ar[dd]^t\\
  \\
  A\ar[rr]_s&&C
}  
\end{aligned}
\end{equation}
For morphisms $\sigma\colon C\to A$ and $\tau\colon C\to B$ of
$\mathcal{C}$, the pushout (if it exists) is denoted by
\begin{equation}
  A\textstyle\pushout{\sigma}{C}{\tau}B:=\lim_{\rightarrow}(\xymatrix@1{A&C\ar[l]_\sigma\ar[r]^\tau&B}),
\end{equation}
\ie\
\begin{equation}
\begin{aligned}
\xymatrix{
  C\ar[rr]^\tau\ar[dd]_\sigma&&B\ar[dd]^{\imath_2}\\
  \\
  A\ar[rr]_{\imath_1}&&A\pushout{\sigma}{C}{\tau}B
}
\end{aligned}
\end{equation}

If for some object $A$ of $\mathcal{C}$, the pullback
$A\pullback{\id_A}{A}{\id_A}A$ exists, then there is a unique
\emph{diagonal morphism} $\delta\colon A\to
A\pullback{\id_A}{A}{\id_A}A$ such that
$\delta\,p_1=\id_A=\delta\,p_2$. A similar result for the product is
probably more familiar: if for some object $A$ of $\mathcal{C}$, the
product $A\prod A$ exists, then there is a unique \emph{diagonal
morphism} $\delta\colon A\to A\prod A$ such that $\delta
p_1=\id_A=\delta p_2$.

Dually, if the pushout $A\pushout{\id_A}{A}{\id_A}A$ exists, there is
a unique \emph{codiagonal} $\delta^\op\colon
A\pushout{\id_A}{A}{\id_A}A\to A$ such that
$\imath_1\,\delta^\op=\id_A=\imath_2\,\delta^\op$, and if the
coproduct $A\coprod A$ exists, there is a unique \emph{codiagonal}
$\delta^\op\colon A\coprod A\to A$ such that $\imath_1\delta^\op =
\id_A = \imath_2\delta^\op$.

Let $A$, $B$, $A^\prime$ and $B^\prime$ be objects of $\mathcal{C}$
such that the products $A\prod B$ and $A^\prime\prod B^\prime$ exist,
and let $f_A\colon A\to A^\prime$ and $f_B\colon B\to B^\prime$ be
morphisms of $\mathcal{C}$. Then there exists a unique morphism
$(f_A;f_B)\colon A\prod B\to A^\prime\prod B^\prime$ such that
$(f_A;f_B)p_1^\prime=p_1f_A$ and $(f_A;f_B)p_2^\prime=p_2f_B$.

For pullbacks, there is the following refinement of this
construction. Let $s\colon A\to C$, $t\colon B\to C$ and
$s^\prime\colon A^\prime\to C^\prime$, $t^\prime\colon B^\prime\to
C^\prime$ be morphisms of $\mathcal{C}$ such that both pullbacks
$A\pullback{s}{C}{t}B$ and
$A^\prime\pullback{s^\prime}{C^\prime}{t^\prime}B$ exist. Let
$f_A\colon A\to A^\prime$, $f_B\colon B\to B^\prime$ and $f_C\colon
C\to C^\prime$ be morphisms such that $f_A\,s^\prime=s\,f_C$ and
$f_B\,t^\prime=t\,f_C$,
\begin{equation}
\begin{aligned}
\xymatrix{
  A\pullback{s}{C}{t}B\ar[rr]^{p_2}\ar[dd]_{p_1}\ar@{.>}[drrr]_(.4){{(f_A;f_B)}_{f_C}}
    &&B\ar[dd]_t\ar[drrr]^{f_B}\\
  &&&A^\prime\pullback{s^\prime}{C^\prime}{t^\prime}B^\prime\ar[rr]_{p_2^\prime}\ar[dd]^{p_1^\prime}
   &&B^\prime\ar[dd]^{t^\prime}\\
  A\ar[rr]^s\ar[drrr]_{f_A}
    &&C\ar[drrr]^(.6){f_C}\\
  &&&A^\prime\ar[rr]_{s^\prime}&&C^\prime
}
\end{aligned}
\end{equation}
Then there exists a unique morphism ${(f_A;f_B)}_{f_C}\colon
A\pullback{s}{C}{t}B\to
A^\prime\pullback{s^\prime}{C^\prime}{t^\prime}B^\prime$ with the
property that ${(f_A;f_B)}_{f_C}\,p_1^\prime = p_1\,f_A$ and
${(f_A;f_B)}_{f_C}\,p_2^\prime = p_2\,f_B$.

Dually, for objects $A$, $B$, $A^\prime$ and $B^\prime$ of
$\mathcal{C}$ for which the coproducts $A\coprod B$ and
$A^\prime\coprod B^\prime$ exist and for morphisms $f_A\colon A\to
A^\prime$ and $f_B\colon B\to B^\prime$, there exists a unique
morphism $[f_A;f_B]\colon A\coprod B\to A^\prime\coprod B^\prime$ such
that $\imath_1[f_A;f_B]=f_A\imath_1^\prime$ and
$\imath_2[f_A;f_B]=f_B\imath_2^\prime$.

For pushouts, we have the following refinement. For morphisms
$\sigma\colon C\to A$, $\tau\colon C\to B$ and $\sigma^\prime\colon
C^\prime\to A^\prime$, $\tau^\prime\colon C^\prime\to B^\prime$ for
which the pushouts $A\pushout{\sigma}{C}{\tau}B$ and
$A^\prime\pushout{\sigma^\prime}{C^\prime}{\tau^\prime}B^\prime$
exist, and morphisms $f_A\colon A\to A^\prime$, $f_B\colon B\to
B^\prime$ and $f_C\colon C\to C^\prime$ that satisfy
$\sigma^\prime\,f_A=f_C\,\sigma$ and $\tau^\prime\,f_B=f_C\,\tau$,
there is a unique morphism ${[f_A;f_B]}_{f_C}\colon
A\pushout{\sigma}{C}{\tau}B\to
A^\prime\pushout{\sigma^\prime}{C^\prime}{\tau^\prime}B^\prime$ such
that $\imath_1\,{[f_A;f_B]}_{f_C}=f_A\,\imath_1^\prime$ and
$\imath_2\,{[f_A;f_B]}_{f_C}=f_B\,\imath_2$.

Notice that for suitable morphisms, the object (if it exists),
\begin{equation}
  A\textstyle\pullback{s}{C}{t}B\textstyle\pullback{u}{D}{v}{E}:=\lim_{\leftarrow}{
  (\xymatrix@1{A\ar[r]^s&C&B\ar[l]_t\ar[r]^u&D&E\ar[l]_v}}),
\end{equation}
is naturally isomorphic to both
$\bigl(A\pullback{s}{C}{t}B\bigr)\pullback{p_2u}{D}{v}E$ and
$A\pullback{s}{C}{p_1t}\bigl(B\pullback{u}{D}{v}E\bigr)$, where in the
first case, $p_2$ is the second projection associated with the
pullback in parentheses, and in the second case, $p_1$ is the first
projection from the pullback in parentheses. Dually, the object,
\begin{equation}
  A\textstyle\pushout{\sigma}{C}{\tau}B\textstyle\pushout{\phi}{D}{\chi}{E}:=\lim_{\rightarrow}{
  (\xymatrix@1{A&C\ar[l]_\sigma\ar[r]^\tau&B&D\ar[l]_\phi\ar[r]^\chi&E}}),
\end{equation}
(if it exists) is naturally isomorphic to both
$\bigl(A\pushout{\sigma}{C}{\tau}B)\pushout{\phi\imath_2}{D}{\chi}E$
and
$A\pushout{\sigma}{C}{\tau\imath_1}\bigl(B\pushout{\phi}{D}{\chi}E\bigr)$.

We use the term `$2$-category' for a strict $2$-category as opposed to
a bicategory. Our terminology for $2$-categories, $2$-functors,
$2$-natural transformations, $2$-equivalences and $2$-adjunctions is
the same as in~\cite{Bo94}.

\subsection{Internal categories}
\label{sect_internalization}

Many $2$-categories that appear in this article, can be constructed by
\emph{internalization}. The concept of internalization goes back to
Ehresmann~\cite{Eh66}. Here we summarize the key definitions and
results. For more details and proofs, see, for example~\cite{Bo94}.

\begin{definition}
\label{def_internalcat}
Let $\mathcal{C}$ be a finitely complete category.
\begin{myenumerate}
\item
  An \emph{internal category} $C=(C_0,C_1,s,t,\imath,\circ)$ in
  $\mathcal{C}$ consists of objects $C_0$ and $C_1$ of $\mathcal{C}$
  with morphisms $s,t\colon C_1\to C_0$, $\imath\colon C_0\to C_1$ and
  $\circ\colon C_1\pullback{s}{C_0}{t}C_1\to C_1$ of $\mathcal{C}$
  such that the following diagrams commute,
\begin{equation}
\label{eq_intcat1}
\begin{aligned}
\xymatrix{
  C_0\ar[rr]^\imath\ar[rrdd]_{\id_{C_0}}&&C_1\ar[dd]^s\\
  \\
  &&C_0
}
\end{aligned}\qquad\qquad
\begin{aligned}
\xymatrix{
  C_0\ar[rr]^\imath\ar[rrdd]_{\id_{C_0}}&&C_1\ar[dd]^t\\
  \\
  &&C_0
}
\end{aligned}
\end{equation}
\begin{equation}
\label{eq_intcat2}
\begin{aligned}
\xymatrix{
  C_1\pullback{s}{C_0}{t}C_1\ar[rr]^\circ\ar[dd]_{p_2}&&C_1\ar[dd]^t\\
  \\
  C_1\ar[rr]_t&&C_0
}
\end{aligned}\qquad\qquad
\begin{aligned}
\xymatrix{
  C_1\pullback{s}{C_0}{t}C_1\ar[rr]^\circ\ar[dd]_{p_1}&&C_1\ar[dd]^s\\
  \\
  C_1\ar[rr]_s&&C_0
}
\end{aligned}
\end{equation}
\begin{equation}
\label{eq_intcat3}
\begin{aligned}
\xymatrix{
  C_1\pullback{s}{C_0}{t}C_1\pullback{s}{C_0}{t}C_1\ar[rr]^{{(\circ;\id_{C_1})}_{\id_{C_0}}}
    \ar[dd]_{{(\id_{C_1};\circ)}_{\id_{C_0}}}&&
    C_1\pullback{s}{C_0}{t}C_1\ar[dd]^\circ\\
  \\
  C_1\pullback{s}{C_0}{t}C_1\ar[rr]_\circ&&C_1
}
\end{aligned}
\end{equation}
\begin{equation}
\label{eq_intcat4}
\begin{aligned}
\xymatrix{
  C_0\pullback{\id_{C_0}}{C_0}{t}C_1\ar[rr]^{{(\imath;\id_{C_1})}_{\id_{C_0}}}\ar[rrdd]_{p_2}&&
    C_1\pullback{s}{C_0}{t}C_1\ar[dd]_\circ&&
    C_1\pullback{s}{C_0}{\id_{C_0}}C_0\ar[ll]_{{(\id_{C_1};\imath)}_{\id_{C_0}}}\ar[ddll]^{p_1}\\
  \\
  && C_1
}
\end{aligned}
\end{equation}
\item
  Let $C=(C_0,C_1,s,t,\imath,\circ)$ and
  $C^\prime=(C_0^\prime,C_1^\prime,s^\prime,t^\prime,\imath^\prime,\circ^\prime)$
  be internal categories in $\mathcal{C}$. An \emph{internal functor}
  $F=(F_0,F_1)\colon C\to C^\prime$ in $\mathcal{C}$ consists of
  morphisms $F_0\colon C_0\to C_0^\prime$ and $F_1\colon C_1\to
  C_1^\prime$ of $\mathcal{C}$ such that the following diagrams
  commute,
\begin{equation}
\label{eq_intfunc1}
\begin{aligned}
\xymatrix{
  C_1\ar[rr]^s\ar[dd]_{F_1}&&C_0\ar[dd]^{F_0}\\
  \\
  C_1^\prime\ar[rr]_{s^\prime}&&C_0^\prime
}
\end{aligned}\qquad\qquad
\begin{aligned}
\xymatrix{
  C_1\ar[rr]^t\ar[dd]_{F_1}&&C_0\ar[dd]^{F_0}\\
  \\
  C_1^\prime\ar[rr]_{t^\prime}&&C_0^\prime
}
\end{aligned}
\end{equation}
\begin{equation}
\label{eq_intfunc2}
\begin{aligned}
\xymatrix{
  C_0\ar[rr]^\imath\ar[dd]_{F_0}&&C_1\ar[dd]^{F_1}\\
  \\
  C_0^\prime\ar[rr]_{\imath^\prime}&&C_1^\prime
}
\end{aligned}\qquad\qquad
\begin{aligned}
\xymatrix{
  C_1\pullback{s}{C_0}{t}C_1\ar[rr]^\circ\ar[dd]_{{(F_1;F_1)}_{F_0}}&&C_1\ar[dd]^{F_1}\\
  \\
  C_1^\prime\pullback{s^\prime}{C_0^\prime}{t^\prime}C_1^\prime\ar[rr]_{\circ^\prime}&&C_1^\prime  
}
\end{aligned}
\end{equation}
\item
  Let $C=(C_0,C_1,s,t,\imath,\circ)$,
  $C^\prime=(C_0^\prime,C_1^\prime,s^\prime,t^\prime,\imath^\prime,\circ^\prime)$
  and $C^{\prime\prime}=(C_0^{\prime\prime},C_1^{\prime\prime},
  s^{\prime\prime},t^{\prime\prime},\imath^{\prime\prime},\circ^{\prime\prime})$
  be internal categories in $\mathcal{C}$ and $F\colon C\to C^\prime$
  and $G\colon C^\prime\to C^{\prime\prime}$ be internal functors. The
  \emph{composition} of $F$ and $G$ is the internal functor $FG\colon
  C\to C^{\prime\prime}$ which is defined by ${(FG)}_0:=F_0G_0$ and
  ${(FG)}_1:=F_1G_1$. 
\item
  Let $C$ be an internal category in $\mathcal{C}$. The \emph{identity
  internal functor} $\id_C\colon C\to C$ is defined by
  ${(\id_C)}_0:=\id_{C_0}$ and ${(\id_C)}_1:=\id_{C_1}$.
\item
  Let $F,\tilde F\colon C\to C^\prime$ be internal functors between
  internal categories $C=(C_0,C_1,s,t,\imath,\circ)$ and
  $C^\prime=(C_0^\prime,C_1^\prime,s^\prime,t^\prime,\imath^\prime,\circ^\prime)$
  in $\mathcal{C}$. An \emph{internal natural transformation}
  $\eta\colon F\Rightarrow\tilde F$ is a morphism $\eta\colon C_0\to
  C_1^\prime$ of $\mathcal{C}$ such that the following diagrams
  commute,
\begin{equation}
\label{eq_intnat1}
\begin{aligned}
\xymatrix{
  C_0\ar[rr]^\eta\ar[rrdd]_{F_0}&&C_1^\prime\ar[dd]^{s^\prime}\\
  \\
  &&C_0^\prime
}
\end{aligned}\qquad\qquad
\begin{aligned}
\xymatrix{
  C_0\ar[rr]^\eta\ar[rrdd]_{\tilde F_0}&&C_1^\prime\ar[dd]^{t^\prime}\\
  \\
  &&C_0^\prime
}
\end{aligned}
\end{equation}
\begin{equation}
\label{eq_intnat2}
\xymatrix{
  C_1\ar[rr]^{\delta{(t\eta;\tilde F_1)}_{t\tilde F_0}}
    \ar[dd]_{\delta{(F_1;s\eta)}_{sF_0}}
    &&C_1^\prime\pullback{s^\prime}{C_0^\prime}{t^\prime}C_1^\prime\ar[dd]^{\circ^\prime}\\
  \\
  C_1^\prime\pullback{s^\prime}{C_0^\prime}{t^\prime}C_1^\prime\ar[rr]_{\circ^\prime}&&C_1^\prime
}  
\end{equation}
\item
  Let $\eta\colon F\Rightarrow\tilde F$ and $\theta\colon\tilde
  F\Rightarrow\hat F$ be internal natural transformations between
  internal functors $F,\tilde F,\hat F\colon C\to C^\prime$. The
  \emph{vertical composition} (or just composition) of $\eta$ and
  $\theta$ is the internal natural transformation
  $\theta\circ\eta\colon F\Rightarrow\hat F$ defined by the following
  composition,
\begin{equation}
\xymatrix@1{
  C_0\ar[r]^(.3){\delta}
    &C_0\pullback{\id_{C_0}}{C_0}{\id_{C_0}}C_0\ar[rr]^(.55){{(\eta;\theta)}_{\tilde F_0}}
    &&C_1^\prime\pullback{s^\prime}{C_0^\prime}{t^\prime}C_1^\prime\ar[r]^(.6){\circ^\prime}
    &C_1^\prime
}
\end{equation}
  Note that the vertical composition `$\circ$' is read from right to
  left\footnote{This is a deviation from~\cite{Pf03} which, however,
  will turn out to be more natural in a future work.}.
\item
  Let $F\colon C\to C^\prime$ be an internal functor in
  $\mathcal{C}$. The \emph{identity internal natural transformation}
  $\id_F\colon F\Rightarrow F$ is defined by $\id_F:=F_0\imath^\prime$
  or by $\imath F_1$.
\item
  Let $\eta\colon F\Rightarrow\tilde F$ and $\tau\colon
  G\Rightarrow\tilde G$ be internal natural transformations between
  internal functors $F,\tilde F\colon C\to C^\prime$ and $G,\tilde
  G\colon C^\prime\to C^{\prime\prime}$. The \emph{horizontal
  composition} (or Godement product) of $\eta$ and $\tau$ is the
  internal natural transformation $\eta\cdot\tau\colon
  FG\Rightarrow\tilde F\tilde G$ defined by the composition,
\begin{equation}
\xymatrix@1{
  C_0\ar[r]^(.25){\delta}
    &C_0\dpullback{\id_{C_0}}{C_0}{\id_{C_0}}C_0\ar[rr]^{{(F_0;\eta)}_{F_0}}
    &&C_0^\prime\dpullback{\id_{C_0^\prime}}{C_0^\prime}{t^\prime}C_1^\prime\ar[rr]^{{(\tau;\tilde G_1)}_{\tilde G_0}}
    &&C_1^{\prime\prime}\dpullback{s^{\prime\prime}}{C_0^{\prime\prime}}{t^{\prime\prime}}C_1^{\prime\prime}
      \ar[r]^(.65){\circ^{\prime\prime}}
    &C_1^{\prime\prime}
}
\end{equation}
or equivalently by,
\begin{equation}
\xymatrix@1{
  C_0\ar[r]^(.25){\delta}
    &C_0\dpullback{\id_{C_0}}{C_0}{\id_{C_0}}C_0\ar[rr]^{{(\eta;\tilde F_0)}_{\tilde F_0}}
    &&C_1^\prime\dpullback{s^\prime}{C_0^\prime}{\id_{C_0^\prime}}C_0^\prime\ar[rr]^{{(G_1;\tau)}_{G_0}}
    &&C_1^{\prime\prime}\dpullback{s^{\prime\prime}}{C_0^{\prime\prime}}{t^{\prime\prime}}C_1^{\prime\prime}
      \ar[r]^(.65){\circ^{\prime\prime}}
    &C_1^{\prime\prime}
}
\end{equation}
\end{myenumerate}
\end{definition}

\begin{theorem}
Let $\mathcal{C}$ be a finitely complete category. There is a
$2$-category $\Cat(\mathcal{C})$ whose objects are internal categories
in $\mathcal{C}$, whose morphisms are internal functors and whose
$2$-morphisms are internal natural transformations.
\end{theorem}

The following example is the motivation for the precise details of the
definition of $\Cat(\mathcal{C})$ presented above.

\begin{example}
The category $\Set$ is finitely complete. $\Cat(\Set)$ is the
$2$-category $\Cat$ of small categories with functors and natural
transformations.
\end{example}

\begin{remark}
For a technical subtlety that arises because limit objects are
specified only up to (unique) isomorphism, we refer to
Remark~\ref{rem_technical} in the Appendix.

For generic finitely complete $\mathcal{C}$, the $1$-category
underlying $\Cat(\mathcal{C})$ is studied in the context of
\emph{essentially algebraic theories}~\cite{Fr72}, going back to the
work of Lawvere on functorial semantics~\cite{La63}. For more details
and references, see, for example~\cite{BaWe83}. The \emph{theory of
categories} $\Th(\Cat)$ is the smallest finitely complete category
that contains objects $C_0$, $C_1$ and morphisms $s,t\colon C_1\to
C_0$, $\imath\colon C_1\to C_0$ and $\circ\colon
C_1\textstyle\pullback{s}{C_0}{t}C_1\to C_1$ such that the
relations~\eqref{eq_intcat1}--\eqref{eq_intcat4} hold. A \emph{model}
of $\Th(\Cat)$ is a finite-limit preserving functor
$F\colon\Th(\Cat)\to\mathcal{C}$ into some finitely complete category
$\mathcal{C}$. If one denotes by
$\Mod(\Th(\Cat),\mathcal{C}):=\fcCat(\Th(\Cat),\mathcal{C})$ the
category of finite-limit preserving (left exact) functors
$\Th(\Cat)\to\mathcal{C}$ with their natural transformations, then
$\Mod(\Th(\Cat),\mathcal{C})\simeq\Cat(\mathcal{C})$ are equivalent as
$1$-categories, in particular $\Mod(\Th(\Cat),\Set)\simeq\Cat$, which
justifies the terminology \emph{theory of categories} for $\Th(\Cat)$.

Usually, these techniques are used in order to study various algebraic
or essentially algebraic theories. One defines, for example,
$\Th(\Grp)$ such that $\Mod(\Th(\Grp),\Set)\simeq\Grp$, \etc. The
framework of internalization is also used in topos theory in order to
replace the category $\Set$ by \emph{more general} finitely complete
categories.

In the following, we are in addition interested in the $2$-categorical
structure of $\Cat(\mathcal{C})$ which cannot be directly seen from
$\Mod(\Th(\Cat),\mathcal{C})$, and we vary the base category
$\mathcal{C}$ over familiar categories of algebraic structures
\emph{more special} than $\Set$, for example, $\mathcal{C}=\Grp$. Then
$\Cat(\mathcal{C})$ forms a $2$-category whose objects turn out to be
(usually the strict versions of) novel higher-dimensional algebraic
structures. Their higher-dimensional nature is directly related to the
fact that $\Cat(\mathcal{C})$ forms a $2$-category and thus exhibits
one more level of structure than the $1$-category $\mathcal{C}$ which
we have initially supplied.
\end{remark}

We can perform the construction $\Cat(\mathcal{C})$ for various
finitely complete base categories $\mathcal{C}$. The following
propositions show what happens if these different base categories are
related by finite-limit preserving functors and by their natural
transformations.

\begin{proposition}
Let $\mathcal{C}$ and $\mathcal{D}$ be finitely complete categories
and $T\colon\mathcal{C}\to\mathcal{D}$ be a functor that preserves
finite limits. Then there is a $2$-functor
$\Cat(T)\colon\Cat(\mathcal{C})\to\Cat(\mathcal{D})$ given as follows.
\begin{myenumerate}
\item
  $\Cat(T)$ associates with each internal category
  $C=(C_0,C_1,s,t,\imath,\circ)$ in $\mathcal{C}$, the internal
  category $\Cat(T)[C]:=(TC_0,TC_1,Ts,Tt,T\imath,T\circ)$ in
  $\mathcal{D}$.
\item
  Let $C$ and $C^\prime$ be internal categories in
  $\mathcal{C}$. $\Cat(T)$ associates with each internal functor
  $F=(F_0,F_1)\colon C\to C^\prime$ in $\mathcal{C}$, the internal
  functor $\Cat(T)[F]\colon\Cat(T)[C]\to\Cat(T)[C^\prime]$ in
  $\mathcal{D}$ which is given by the following morphisms of
  $\mathcal{D}$,
\begin{eqnarray}
  {(\Cat(T)[F])}_0&=&TF_0\colon TC_0\to TC_0^\prime,\\
  {(\Cat(T)[F])}_1&=&TF_1\colon TC_1\to TC_1^\prime.
\end{eqnarray}
\item
  Let $C$ and $C^\prime$ be internal categories and $F,F^\prime\colon
  C\to C^\prime$ be internal functors in $\mathcal{C}$. $\Cat(T)$
  associates with each internal natural transformation $\eta\colon
  F\Rightarrow F^\prime$, the internal natural transformation
  $\Cat(T)[\eta]\colon\Cat(T)[F] \Rightarrow\Cat(T)[F^\prime]$ in
  $\mathcal{D}$ which is defined by the following morphism of
  $\mathcal{D}$,
\begin{equation}
  \Cat(T)[\eta]=T\eta\colon TC_0\to TC_1^\prime.
\end{equation}
\end{myenumerate}
\end{proposition}

\begin{proposition}
Let $\mathcal{C}$ and $\mathcal{D}$ be finitely complete categories,
$T,\tilde T\colon\mathcal{C}\to\mathcal{D}$ be functors that preserve
finite limits and $\alpha\colon T\Rightarrow\tilde T$ a natural
transformation. Then there is a $2$-natural transformation
$\Cat(\alpha)\colon\Cat(T)\Rightarrow\Cat(\tilde T)$ given as follows.

$\Cat(\eta)$ associates with each internal category
$C=(C_0,C_1,s,t,\imath,\circ)$ in $\mathcal{C}$ the internal functor
${\Cat(\alpha)}_C\colon\Cat(T)[C]\to\Cat(\tilde T)[C]$ in
$\mathcal{D}$ which is given by the following morphisms of
$\mathcal{D}$,
\begin{eqnarray}
  {({\Cat(\alpha)}_C)}_0&=&\alpha_{C_0}\colon TC_0\to\tilde TC_0,\\
  {({\Cat(\alpha)}_C)}_1&=&\alpha_{C_1}\colon TC_1\to\tilde TC_1.
\end{eqnarray}
\end{proposition}

\begin{theorem}
\label{thm_internalize}
Let $\fcCat$ denote the $2$-category of finitely complete small
categories, finite-limit preserving functors and natural
transformations. Let $\TCat$ be the $2$-category of small
$2$-categories, $2$-functors and $2$-natural transformations. Then
$\Cat(-)$ forms a $2$-functor\footnote{It is known~\cite{BaWe83} that
$\Cat(-)$ is a functor $\fcCat\to\fcCat$, but we do note make use of
this fact here.},
\begin{equation}
  \Cat(-)\colon\fcCat\to\TCat.
\end{equation}
\end{theorem}

\begin{corollary}
\label{corr_equivalence}
Let $\mathcal{C}\simeq\mathcal{D}$ be an equivalence of finitely
complete categories provided by the functors
$F\colon\mathcal{C}\to\mathcal{D}$ and
$G\colon\mathcal{D}\to\mathcal{C}$ with natural isomorphisms
$\eta\colon 1_\mathcal{C}\Rightarrow FG$ and $\epsilon\colon
GF\Rightarrow 1_\mathcal{D}$. Then there is a $2$-equivalence of the
$2$-categories $\Cat(\mathcal{C})\simeq\Cat(\mathcal{D})$ given by the
$2$-functors $\Cat(F)\colon\Cat(\mathcal{C})\to\Cat(\mathcal{D})$ and
$\Cat(G)\colon\Cat(\mathcal{D})\to\Cat(\mathcal{C})$ with the
$2$-natural isomorphisms $\Cat(\eta)\colon
1_{\Cat(\mathcal{C})}\Rightarrow\Cat(F)\Cat(G)$ and
$\Cat(\epsilon)\colon\Cat(G)\Cat(F)\Rightarrow 1_{\Cat(\mathcal{D})}$.
\end{corollary}

We can now consider any diagram in $\fcCat$, in particular any diagram
involving finitely complete categories of familiar algebraic
structures, finite-limit preserving functors and their natural
transformations. Theorem~\ref{thm_internalize} guarantees that the
$\Cat(-)$-image of any such diagram is a valid diagram of
$2$-categories, $2$-functors and $2$-natural transformations. This
idea is the key to the present article and provides us with the
desired $2$-functors between our $2$-categories of higher-dimensional
algebraic structures.

The study of contravariant functors, \ie\ functors that can be written
covariantly as $T\colon\mathcal{D}\to\mathcal{C}^\op$, leads to the
following concept dual to the notion of internal categories.

\begin{definition}
\label{def_cocat}
Let $\mathcal{C}$ be a finitely cocomplete category. An \emph{internal
cocategory} in $\mathcal{C}$ is an internal category in
$\mathcal{C}^\op$. \emph{Internal cofunctors} and \emph{internal
conatural transformations} are internal functors and internal natural
transformations in $\mathcal{C}^\op$. We denote by
$\CoCat(\mathcal{C}):=\Cat(\mathcal{C}^\op)$ the $2$-category of
internal cocategories, cofunctors and conatural transformations in
$\mathcal{C}$.
\end{definition}

More explicitly, an internal cocategory
$C=(C_0,C_1,\sigma,\tau,\Uepsilon,\UDelta)$ in $\mathcal{C}$ consists
of objects $C_0$ and $C_1$ of $\mathcal{C}$ and morphisms
$\sigma,\tau\colon C_0\to C_1$, $\Uepsilon\colon C_1\to C_0$ and
$\UDelta\colon C_1\to C_1\pushout{\sigma}{C_0}{\tau}C_1$ of
$\mathcal{C}$ such that the diagrams dual
to~\eqref{eq_intcat1}--\eqref{eq_intcat4} commute if all the morphisms
are relabelled as follows: $s\mapsto\sigma$; $t\mapsto\tau$;
$\imath\mapsto\Uepsilon$ and $\circ\mapsto\UDelta$. Notice that all
arrows involved in the universal constructions have to be reversed,
too, for example, pullbacks have to be replaced by pushouts and
diagonal morphisms by codiagonal ones.

Given internal cocategories $C$ and $C^\prime$ in $\mathcal{C}$, an
internal cofunctor $F=(F_0,F_1)\colon C\to C^\prime$ consists of
morphisms $F_0\colon C_0^\prime\to C_0$ and $F_1\colon C_1^\prime\to
C_1$ of $\mathcal{C}$ such that the diagrams dual
to~\eqref{eq_intfunc1} and~\eqref{eq_intfunc2} commute.

Given internal cofunctors $F,\tilde F\colon C\to C^\prime$ in
$\mathcal{C}$, an internal conatural transformation is a morphism
$\eta\colon C_1^\prime\to C_0$ such that the diagrams dual
to~\eqref{eq_intnat1} and~\eqref{eq_intnat2} commute.

Notice the counter-intuitive direction of the morphisms $F_0$, $F_1$
and $\eta$ which is a consequence of our definition of $\CoCat(-)$. We
finally remark that we have never called $\Cat(-)$ a
\emph{categorification} since it is not clear in which sense it can be
reversed and whether this would correspond to a form of
decategorification~\cite{BaDo98}.

\subsection{Strict $2$-groups}

Strict $2$-groups form one of the simplest examples of
higher-dimensional algebraic structures. Just as a group can be viewed
as a groupoid with one object, every strict $2$-group gives rise to a
$2$-groupoid with one object. $2$-groups have appeared in the
literature in various contexts and under different names (categorical
group, gr-category, $\mathrm{cat}^1$-group, \etc). For an overview and
a comprehensive list of references, we refer to~\cite{BaLa04}. Strict
$2$-groups can be defined in several different ways. In the present
article, we define a strict $2$-group as an internal category in the
category of groups so that the techniques of
Section~\ref{sect_internalization} are available. Alternatively,
strict $2$-groups are group objects in the category of small
categories, see, for example~\cite{Fo02}.

In order to make the presentation self-contained, let us first recall
the construction of finite limits in the finitely complete category
$\Grp$ of groups. The terminal object is the trivial group $\{e\}$
with the trivial group homomorphisms $G\to\{e\}$; the binary product
of groups $G$ and $H$ is the direct product, $G\prod H=G\times H$,
with the projections $p_1\colon G\times H\to G, (g,h)\mapsto g$ and
$p_2\colon G\times H\to H, (g,h)\mapsto h$; and for group
homomorphisms $f_1,f_2\colon G\to H$, their equalizer is a subgroup of
$G$,
\begin{equation}
\label{eq_equalizergrp}
  \eq(f_1,f_2)=\{\,g\in G\colon\quad f_1(g)=f_2(g)\,\}\subseteq G,
\end{equation}
with its inclusion $e={(\id_G)}|_{\eq(f_1,f_2)}\colon\eq(f_1,f_2)\to
G$.

For group homomorphisms $t\colon G\to K$ and $t\colon H\to K$, we
therefore obtain the pullback,
\begin{equation}
  G\textstyle\pullback{s}{K}{t}H=G\times_K H
  =\{\,(g,h)\in G\times H\colon\quad s(g)=t(h)\,\}\subseteq G\times H,
\end{equation}
with the projections $p_1\colon G\times_KH\to G$, $(g,h)\mapsto g$ and
$p_2\colon G\times_KH\to H$, $(g,h)\mapsto h$.

\begin{definition}
The objects, morphisms and $2$-morphisms of $\TGrp:=\Cat(\Grp)$ are
called \emph{strict $2$-groups}, \emph{homomorphisms} and
\emph{$2$-homomorphisms} of strict $2$-groups, respectively. The
objects of $\fTGrp:=\Cat(\fGrp)$ are called \emph{strict finite
$2$-groups}.
\end{definition}

Examples of strict $2$-groups, their homomorphisms and
$2$-homomorphisms can be constructed from Whitehead's crossed modules
as follows~\cite{Po87}.

\begin{definition}
\begin{myenumerate}
\item
  A \emph{crossed module} $(G,H,\rhd,\del)$ consists of groups $G$ and
  $H$ and group homomorphisms $\del\colon H\to G$ and $G\to\Aut
  H,g\mapsto(h\mapsto g\rhd h)$ that satisfy for all $g\in G$,
  $h,h^\prime\in H$,
\begin{eqnarray}
  \del(g\rhd h)&=&g\del(h)g^{-1},\\
  \del(h)\rhd h^\prime &=& hh^\prime h^{-1}.
\end{eqnarray}
\item
 A \emph{homomorphism}
 $F=(F_G,F_H)\colon(G,H,\rhd,\del)\to(G^\prime,H^\prime,\rhd^\prime,\del^\prime)$
 of crossed modules consists of group homomorphisms $F_G\colon G\to
 G^\prime$ and $F_H\colon H\to H^\prime$ such that,
\begin{equation}
\begin{aligned}
\xymatrix{
  H\ar[rr]^\del\ar[dd]_{F_H}&&G\ar[dd]^{F_G}\\
  \\
  H^\prime\ar[rr]_{\del^\prime}&&G^\prime
}
\end{aligned}
\end{equation}
  commutes and such that for all $g\in G$, $h\in H$,
\begin{equation}
  F_H(g\rhd h)=F_G(g)\rhd^\prime F_H(h).
\end{equation}
\item
  Let $F,\tilde
  F\colon(G,H,\rhd,\del)\to(G^\prime,H^\prime,\rhd^\prime,\del^\prime)$
  be homomorphisms of crossed modules. A \emph{$2$-homomorphism}
  $\eta\colon F\Rightarrow \tilde F$ is a map $\eta_H\colon G\to
  H^\prime$ such that for all $g,g_1,g_2\in G$, $h\in H$,
\begin{eqnarray}
  \eta_H(e)&=&e,\\
  \eta_H(g_1g_2)&=&\eta_H(g_1)(F_G(g_1)\rhd^\prime\eta_H(g_2)),\\
  \del^\prime(\eta(g))&=&\tilde F_G(g){F_G(g)}^{-1},\\
  \eta_H(\del(h))&=&\tilde F_H(h){F_H(h)}^{-1}.
\end{eqnarray}
\end{myenumerate}
\end{definition}

\begin{example}
\begin{myenumerate}
\item
  Let $H$ be a finite group and $G=\Aut H$ its group of
  automorphisms. Choose $\del\colon H\to G,h\mapsto (h^\prime\mapsto
  hh^\prime h^{-1})$, and $g\rhd h:=g(h)$ for $g\in\Aut H$ and $h\in
  H$. Then $(G,H,\rhd,\del)$ forms a crossed module.
\item
  Let $G$ be a finite group and $(V,\rho)$ be a finite-dimensional
  representation of $G$, \ie\ $V$ is a $\K$-vector space and
  $\rho\colon G\to \GL_\K(V)$ a homomorphism of groups. Choose
  $H:=(V,+,0)$ to be the additive group underlying the vector space,
  $g\rhd h:=\rho(g)[h]$ for $g\in G$, $h\in H$, and $\del\colon H\to
  G, h\mapsto e$. Then $(G,H,\rhd,\del)$ forms a crossed module.
\item
  More examples are given in~\cite{BaLa04,Fo03}. 
\end{myenumerate}
\end{example}

\begin{theorem}
\label{thm_xmod}
\begin{myenumerate}
\item
  There is a $2$-category $\XMod$ whose objects are crossed modules,
  whose morphisms are homomorphisms of crossed modules and whose
  $2$-morphisms are $2$-homomorphisms of crossed modules.
\item
  The $2$-categories $\TGrp$ and $\XMod$ are $2$-equivalent.
\end{myenumerate}
\end{theorem}

In order to obtain examples of strict $2$-groups from crossed modules,
we need the explicit form of one of the $2$-functors involved in this
$2$-equivalence, $T\colon\XMod\to\TGrp$.
\begin{myenumerate}
\item
  $T$ associates with each crossed module $(G,H,\rhd,\del)$ the strict
  $2$-group $(G_0,G_1,s,t,\imath,\circ)$ defined as follows. The
  groups are $G_0:=G$ and $G_1:=H\rtimes G$ where the semidirect
  product uses the multiplication $(h_1,g_1)\cdot
  (h_2,g_2):=(h_1(g_1\rhd h_2),g_1g_2)$. The group homomorphisms are
  given by $s\colon H\rtimes G\to G, (h,g)\mapsto g$; $t\colon
  H\rtimes G\to G, (h,g)\mapsto\del(h)g$; $\imath\colon G\to H\rtimes
  G, g\mapsto (e,g)$ and $\circ\colon(H\rtimes G)\times(H\rtimes G)\to
  H\rtimes G, ((h_1,g_1),(h_2,g_2))\mapsto
  (h_1,g_1)\circ(h_2,g_2):=(h_1h_2,g_1)$, defined whenever
  $g_1=s(h_1,g_1)=t(h_2,g_2)=\del(h_2)g_2$.
\item
  For each homomorphism $F=(F_G,F_H)\colon
  (G,H,\rhd,\del)\to(G^\prime,H^\prime,\rhd^\prime,\del^\prime)$ of
  crossed modules, there is a homomorphism of strict $2$-groups
  $TF=(F_0,F_1)$ given by $F_0:=F_G\colon G_0\to G_0^\prime$ and
  $F_1:=F_H\times F_G\colon H\rtimes G\to H^\prime\rtimes G^\prime$.
\item
  Let $F,\tilde F\colon
  (G,H,\rhd,\del)\to(G^\prime,H^\prime,\rhd^\prime,\del^\prime)$ be
  homomorphisms of crossed modules. For each $2$-homomorphism
  $\eta\colon F\Rightarrow \tilde F$, there is a $2$-homomorphism
  $T\eta\colon TF\Rightarrow T\tilde F$ of strict $2$-groups given by
  $T\eta\colon G_0\to G_1^\prime=H^\prime\rtimes G^\prime, g\mapsto
  (\eta_H(g),F_G(g))$.
\end{myenumerate}

From the proof of Theorem~\ref{thm_xmod} (see, for
example~\cite{Fo02,Fo03}), one sees that in any internal category
$(G_0,G_1,s,t,\imath,\circ)$ in $\Grp$, the group homomorphism
$\circ\colon G_1\pullback{s}{G_0}{t}G_1\to G_1$ is already uniquely
determined by $s,t,\imath$ and by the group structures of $G_0$ and $G_1$
and, moreover, that each element $g\in G_1$ has a unique inverse
$g^\times\in G_1$ with respect to `$\circ$'. In other words, an
internal category in $\Grp$ is actually an internal groupoid in
$\Grp$.

This result holds not only for internal categories in $\Grp$, but for
internal categories in the category $\Grp(\mathcal{C})$ of group
objects in any category $\mathcal{C}$ that has all finite products
(\cf~Appendix~\ref{app_gpobj}). This result can be stated as follows.

\begin{proposition}
\label{prop_vertical}
Let $\mathcal{C}$ be a category with finite products and
$(G_0,G_1,s,t,\imath,\circ)$ be an internal category\footnote{This
includes the assumption that the required pullbacks exist in
$\Grp(\mathcal{C})$. A sufficient condition is that $\mathcal{C}$ is
finitely complete.} in $\Grp(\mathcal{C})$.
\begin{myenumerate}
\item
  The morphism $\circ\colon G_1\pullback{s}{G_0}{t}G_1\to G_1$ is of
  the form
\begin{equation}
\label{eq_verticalunique}
  g\circ \tilde g=g{(\imath(s(g)))}^{-1}\tilde g,
\end{equation}
  for all $g,\tilde g\in G_1$ that satisfy $s(g)=t(\tilde
  g)$. Conversely, given only the data $(G_0,G_1,s,t,\imath)$,
  equation~\eqref{eq_verticalunique} defines a morphism $\circ\colon
  G_1\pullback{s}{G_0}{t}G_1\to G_1$ that satisfies~\eqref{eq_intcat2}
  to~\eqref{eq_intcat4}.

  In order to simplify the notation in~\eqref{eq_verticalunique}, we
  have pretended that $\mathcal{C}$ is a subcategory of $\Set$ and so
  we can write down this equation for elements.
\item
  There is a morphism $\xi\colon G_1\to G_1$, $g\mapsto g^\times$ such
  that
\begin{eqnarray}
  t(g^\times)&=&s(g),\\
  s(g^\times)&=&t(g),\\
  g^\times\circ g&=&\imath(s(g)),\\
  g\circ g^\times&=&\imath(t(g)).
\end{eqnarray}
  It is given by $\xi(g)=\imath(s(g))g^{-1}\imath(t(g))$ and satisfies
\begin{equation}
  \xi(g\circ \tilde g) = \xi(\tilde g)\circ \xi(g),
\end{equation}
for all $g,\tilde g\in G_1$ for which $s(g)=t(\tilde g)$. Again, we
have pretended that $\mathcal{C}$ is a subcategory of $\Set$.
\end{myenumerate}
\end{proposition}

The map $\xi\colon G_1\to G_1$ which associates with each element
$g\in G_1$ its inverse $g^\times$ with respect to the vertical
composition `$\circ$' and the contravariant counterpart of this map
(Remark~\ref{rem_vertical2} below), are responsible for the functorial
antipode in the Hopf category of representations of the $2$-group
$(G_0,G_1,s,t,\imath,\circ)$. This is shown in
Proposition~\ref{prop_antipode} and further illustrated in
Section~\ref{sect_individual} below.

%
\section{Cocommutative trialgebras}
%
\label{sect_coctrialg}

Given some group $G$, its group algebra $\K[G]$ forms a cocommutative
Hopf algebra. In this section, we use the technique of internalization
in order to construct the analogue of the group algebra for strict
$2$-groups. This gives rise to a novel higher-dimensional algebraic
structure which is defined as an internal category in the category of
cocommutative Hopf algebras. We call this a \emph{cocommutative
trialgebra} for reasons that are explained below. For general
background on coalgebras, bialgebras and Hopf algebras, we refer
to~\cite{Ma95}.

\begin{definition}
The functor $\K[-]\colon\Grp\to\cocHopf_\K$ is defined as follows. It
associates with each group $G$ its group algebra $\K[G]$. This is the
free vector space over the set $G$ equipped with the structure of a
cocommutative Hopf algebra $(\K[G],\mu,\eta,\Delta,\epsilon,S)$ using
the multiplication $\mu\colon\K[G]\otimes\K[G]\to\K[G]$, defined on
basis elements $g,h\in G$ by $g\otimes h\mapsto gh$ (group
multiplication), the unit $\eta\colon\K\to\K[G]$, $1\mapsto e$ (group
unit), comultiplication $\Delta\colon\K[G]\to\K[G]\otimes\K[G]$,
$g\mapsto g\otimes g$ (group-like), counit $\epsilon\colon\K[G]\to\K$,
$g\mapsto 1$, and antipode $S\colon\K[G]\to\K[G]$, $g\mapsto g^{-1}$
(group inverse). The functor $\K[-]$ associates with each group
homomorphism $f\colon G\to H$ the bialgebra homomorphism
$\K[f]\colon\K[G]\to\K[H]$ which is the $\K$-linear extension of $f$.
\end{definition}

Given some strict $2$-group $(G_0,G_1,s,t,\imath,\circ)$, the idea is
to apply $\K[-]$ to $G_0$ and $G_1$ and to all maps
$s,t,\imath,\circ$. The result is a structure $(H_0,H_1,\hat s,\hat
t,\hat\imath,\hat\circ)$ consisting of cocommutative Hopf algebras
$H_0$ and $H_1$ with various bialgebra homomorphisms. In the
following, we show that such a structure can alternatively be defined
as an internal category in the category of cocommutative Hopf
algebras. Definition~\ref{def_internalcat} refers to several universal
constructions such as pullbacks and diagonal morphisms which are all
constructed from finite limits. It is therefore sufficient to show
that $\K[-]$ preserves all finite limits
(\cf~Theorem~\ref{thm_internalize}).

\subsection{Finite limits in the category of cocommutative Hopf algebras}
\label{sect_cochopflimit}

We first recall the construction of finite limits in the category
$\cocHopf_\K$. The proofs of the following propositions are
elementary.

\begin{proposition}
The terminal object in any of the categories $\Hopf_\K$,
$\cocHopf_\K$, $\cocBi_\K$ and $\cocCo_\K$ is given by $\K$
itself. For each object $C$, the unique morphism $C\to\K$ is the
counit operation $\epsilon_C\colon C\to\K$.
\end{proposition}

\begin{proposition}
\label{prop_productcochopf}
The binary product in the category $\cocHopf_\K$ is the tensor product
of Hopf algebras, \ie\ for cocommutative Hopf algebras $A$ and $B$,
\begin{equation}
  A\textstyle\prod B=A\otimes B,
\end{equation}
with the Hopf algebra homomorphisms $p_1=\id_A\otimes\epsilon_B\colon
A\otimes B\to A$ and $p_2=\epsilon_A\otimes\id_B\colon A\otimes B\to
B$. For convenience, we have suppressed the isomorphisms
$A\otimes\K\cong A$ and $\K\otimes B\cong B$.

For each cocommutative Hopf algebra $D$ with homomorphisms $f_1\colon
D\to A$ and $f_2\colon D\to B$, there is a unique Hopf algebra
homomorphism $\phi\colon D\to A\otimes B$ such that $\phi p_1=f_1$ and
$\phi p_2=f_2$. It is given by $\phi:=\Delta_D(f_1\otimes f_2)$ (first
$\Delta_D$, then $f_1\otimes f_2$).
\end{proposition}

\begin{remark}
The binary product in the categories $\cocBi_\K$ and $\cocCo_\K$ is
given precisely by the same construction: the tensor product of
bialgebras or coalgebras, respectively. In the category $\Hopf_\K$ of
all (not necessarily cocommutative) Hopf algebras, the product is in
general \emph{not} the tensor product of Hopf
algebras. Cocommutativity is essential in
Proposition~\ref{prop_productcochopf} in order to show that
$\phi\colon D\to A\otimes B$ is indeed a homomorphism of coalgebras.
\end{remark}

\begin{proposition}
The equalizer of a pair of morphisms (binary equalizer) in the
category $\Hopf_\K$ is given by the coalgebra cogenerated by the
equalizer in $\Vect_\K$, \ie\ for Hopf algebras $D$, $C$ with
homomorphisms $f_1,f_2\colon D\to C$, the equalizer of $f_1$ and $f_2$
is the largest subcoalgebra $\eq(f_1,f_2)$ of $D$ that is contained in
the linear subspace $\ker(f_1-f_2)\subseteq D$, with its inclusion
$e:={(\id_D)}|_{\eq(f_1,f_2)}\colon\eq(f_1,f_2)\to D$. The coalgebra
$\eq(f_1,f_2)$ turns out to be a sub-Hopf algebra of $D$ and $e$ a
homomorphism of Hopf algebras.

For each Hopf algebra $H$ with a homomorphism $\psi\colon H\to D$ such
that $\psi f_1=\psi f_2$, the image $\psi(H)$ forms a coalgebra which
is contained in $\ker(f_1-f_2)$ so that
$\psi(H)\subseteq\eq(f_1,f_2)$. Then there is a unique Hopf algebra
homomorphism $\phi\colon H\to\eq(f_1,f_2)$ such that $\psi=\phi e$. It
is obtained by simply restricting the codomain of $\psi$ to
$\eq(f_1,f_2)$.
\end{proposition}

\begin{remark}
The equalizer in $\cocHopf_\K$, $\cocBi_\K$ and $\cocCo_\K$ is given
precisely by the same construction.
\end{remark}

\begin{corollary}
\label{corr_pullbackcochopf}
The pullback in the category $\cocHopf_\K$ is constructed as
follows. Let $A$, $B$, $C$ be cocommutative Hopf algebras with
homomorphisms $s\colon A\to C$ and $t\colon B\to C$. The pullback
$A\pullback{s}{C}{t}B$ is the largest subcoalgebra of $A\otimes B$
that is contained in the linear subspace $\ker\Phi\subseteq A\otimes
B$. Here $\Phi$ denotes the linear map
$\Phi=(s\otimes\epsilon_B-\epsilon_A\otimes t)\colon A\otimes B\to
C$. The limiting cone is given by the Hopf algebra homomorphisms
$p_1={(\id_A\otimes\epsilon_B)}|_{A\pullback{s}{C}{t}B}\colon
A\pullback{s}{C}{t}B\to A$ and
$p_2={(\epsilon_A\otimes\id_B)}|_{A\pullback{s}{C}{t}B}\colon
A\pullback{s}{C}{t}B\to B$.

If $H$ is a cocommutative Hopf algebra with homomorphisms $f_1\colon
H\to A$ and $f_2\colon H\to B$ such that $f_1s=f_2t$, then there is a
unique homomorphism $\phi\colon H\to A\pullback{s}{C}{t}B$ such that
$\phi p_1=f_1$ and $\phi p_2=f_2$. It is given by restricting the
codomain of $\phi:=\Delta_H(f_1\otimes f_2)\colon H\to A\otimes B$ to
$A\pullback{s}{C}{t}B$. This restriction is possible because $\phi(H)$
is a coalgebra that is contained in the linear subspace $\ker\Phi$
which implies that $\phi(H)\subseteq A\pullback{s}{C}{t}B$.
\end{corollary}

\subsection{Finite-limit preservation of the group algebra functor}

\begin{theorem}
\label{thm_grpalglimit}
The group algebra functor $\K[-]\colon\Grp\to\cocHopf_\K$ preserves
finite limits.
\end{theorem}

\begin{proof}
We verify in a direct calculation that $\K[-]$ preserves the terminal
object, binary products and binary equalizers.
\begin{myenumerate}
\item
  For the terminal object, we note that $\K[\{e\}]\cong\K$ are
  isomorphic as Hopf algebras. 
\item 
  In order to see that $\K[-]$ preserves binary products, consider
  groups $G$ and $H$. We recall that there is an isomorphism of
  Hopf algebras,
\begin{equation}
  \K[G\times H]\to \K[G]\otimes\K[H],\quad (g,h)\mapsto g\otimes h,
\end{equation}
  which shows that $\K[-]$ maps the product object $G\prod H=G\times
  H$ in $\Grp$ to the product object
  $\K[G]\prod\K[H]=\K[G]\otimes\K[H]$ in $\cocHopf_\K$. Let $p_1\colon
  G\times H\to G, (g,h)\mapsto g$ be the first projection, then
  $\K[p_1]\colon\K[G]\otimes\K[H]\to\K[G], g\otimes h\mapsto g$ is
  precisely the map $\K[p_1]=\id_{\K[G]}\otimes\epsilon_{\K[H]}$. An
  analogous result holds for $p_2\colon G\times H\to H, (g,h)\mapsto
  h$. Therefore, $K[-]$ maps the limiting cone in $\Grp$ to the
  limiting cone in $\cocHopf_\K$.
\item
  For equalizers, consider group homomorphisms $f_1,f_2\colon G\to H$
  and denote their equalizer in $\Grp$ by $\eq(f_1,f_2)\subseteq
  G$. We show that,
\begin{equation}
  \K[\eq(f_1,f_2)] = \ker(\K[f_1]-\K[f_2])\subseteq\K[G],
\end{equation}
  \ie\ the kernel of the difference which is the equalizer in
  $\Vect_\K$, already forms a coalgebra, namely $\K[\eq(f_1,f_2)]$
  itself. In order to see this, choose the standard basis $G$ of
  $\K[G]$. Then $g\in \ker(\K[f_1]-\K[f_2])$ if and only if
  $f_1(g)-f_2(g)=0$ in $\K[H]$ which holds if and only if
  $g\in\eq(f_1,f_2)$.
\end{myenumerate}
\end{proof}

\begin{corollary}
The group algebra functor $\K[-]\colon\Grp\to\cocHopf_\K$ preserves
pullbacks. In particular, for groups $G$, $H$ and $K$ and group
homomorphisms $s\colon G\to K$ and $t\colon H\to K$, there is an
isomorphism of Hopf algebras,
\begin{equation}
  \K[G\times_K H] \cong \K[G]\textstyle\pullback{\K[s]}{\K[K]}{\K[t]}\K[H] 
   = \ker\Phi\subseteq\K[G]\otimes\K[H],
\end{equation}
where
$\Phi=(\K[s]\otimes\epsilon_{\K[H]}-\epsilon_{\K[G]}\otimes\K[t])
\colon\K[G]\otimes\K[H]\to\K[K]$.
\end{corollary}

\begin{remark}
If the group algebra functor is viewed as a functor
$\K[-]\colon\Grp\to\Hopf_k$ into the category of all (not necessarily
cocommutative) Hopf algebras, it does not preserve all pullbacks. In
fact, it does not even preserve all binary products. This can be
blamed on the inclusion functor $\cocHopf_\K\to\Hopf_\K$ which does
not preserve all binary products.
\end{remark}

\subsection{Definition of cocommutative trialgebras}
\label{sect_coctrialgdef}

\begin{definition}
The objects, morphisms and $2$-morphisms of the $2$-category,
\begin{equation}
  \cocTriAlg_\K:=\Cat(\cocHopf_\K),
\end{equation} 
are called \emph{strict cocommutative trialgebras}, their
homomorphisms and $2$-homomorphisms, respectively.
\end{definition}

\begin{proposition}
The functor $\K[-]\colon\Grp\to\cocHopf_\K$ gives rise to a
$2$-functor,
\begin{equation}
  \Cat(\K[-])\colon\TGrp\to\cocTriAlg_\K.
\end{equation}
\end{proposition}

This is the main consequence of the preceding section. In particular,
we can use this $2$-functor in order to obtain examples of strict
cocommutative trialgebras from strict $2$-groups and thereby from
Whitehead's crossed modules.

%
%

Let us finally unfold the definition of a strict cocommutative
trialgebra in more detail.

\begin{proposition}
A strict cocommutative trialgebra $H=(H_0,H_1,s,t,\imath,\circ)$
consists of cocommutative Hopf algebras $H_0$ and $H_1$ over $\K$ and
bialgebra homomorphisms $s\colon H_1\to H_0$, $t\colon H_1\to H_0$,
$\imath\colon H_0\to H_1$ and $\circ\colon H_1\pullback{s}{H_0}{t}
H_1\to H_1$ such that~\eqref{eq_intcat1}--\eqref{eq_intcat4} hold
(with $C_j$ replaced by $H_j$).
\end{proposition}

\begin{remark}
\label{rem_trialgebra}
The terminology \emph{trialgebra} for such an internal category $H$ in
$\cocHopf_\K$ originates from the observation that the vector space
$H_1$ is equipped with three linear operations. There is a
multiplication and a comultiplication since $H_1$ forms a Hopf
algebra. In addition, there is another, \emph{partially defined},
multiplication $\circ\colon H_1\pullback{s}{H_0}{t}H_1\to H_1$. From
the construction of the pullback
(Corollary~\ref{corr_pullbackcochopf}), we see that for $h,h^\prime\in
H_1$, the multiplication $h\circ h^\prime$ is defined only if
$h\otimes h^\prime$ lies in the largest subcoalgebra of $H_1\otimes
H_1$ that is contained in the linear subspace,
\begin{equation}
  \ker(s\otimes\epsilon_{H_1}-\epsilon_{H_1}\otimes t)\subseteq H_1\otimes H_1.
\end{equation}
Whenever the multiplication `$\circ$' is defined, \eqref{eq_intcat3} implies
that it is associative.

The purpose of the Hopf algebra $H_0$ and of the homomorphisms
$s,t\colon H_1\to H_0$ which feature in the pullback
$H_1\pullback{s}{H_0}{t}H_1$, is to keep track of precisely when the
multiplication `$\circ$' is defined. This is completely analogous to
the partially defined multiplication $\circ\colon
G_1\times_{G_0}G_1\to G_1$ in a strict $2$-group
$(G_0,G_1,s,t,\imath,\circ)$ where a pair of elements $(g,g^\prime)\in
G_1\times G_1$ is `$\circ$'-composable if and only if the source of
$g$ agrees with the target of $g^\prime$, \ie\ $s(g)=t(g^\prime)$.

In the cocommutative trialgebra $H$, the partially defined
multiplication `$\circ$' has got \emph{local units}. The homomorphism
$\imath\colon H_0\to H_1$ yields the units, and~\eqref{eq_intcat4}
implies the following unit law: the element $\imath(t(h))$ is (up to a
factor) a left-unit for $h\in H_1$, \ie\ $\circ(\imath(t(h))\otimes
h)=\epsilon(h)\,h$, and $\imath(s(h))$ is a right-unit, \ie\
$\circ(h\otimes\imath(s(h)))=h\,\epsilon(h)$. Note that this sort of
units can depend on $s(h)$ and $t(h)$, just as the identity morphisms
in a small category which form the left- or right-units of a given
morphism, depend on the source and target object of that morphism.

All three operations on $H_1$ are compatible in the following way. The
multiplication and comultiplication of $H_1$ are compatible with each
other because $H_1$ forms a bialgebra. The partially defined
multiplication `$\circ$' and the local unit map `$\imath$' are both
homomorphisms of bialgebras which expresses their compatibility with
the Hopf algebra operations of $H_0$ and $H_1$.

In particular, the fact that `$\circ$' forms a homomorphism of
algebras, implies that both multiplications, the globally defined
multiplication `$\cdot$' of the algebra $H_1$ and the partially
defined multiplication `$\circ$' are compatible and satisfy an
interchange law,
\begin{equation}
  \alignidx{(h_1\cdot h_2)\circ (h_1^\prime\cdot h_2^\prime)}
  = \alignidx{(h_1\circ h_1^\prime)\cdot(h_2\circ h_2^\prime)},
\end{equation}
for $\alignidx{h_1,h_2,h_1^\prime,h_2^\prime}\in H_1$, whenever the
partial multiplication `$\circ$' is defined. A possible
Eckmann--Hilton argument~\cite{EcHi62} which would render both
multiplications commutative and equal, is sidestepped by the same
mechanism as in strict $2$-groups: both multiplications can in general
have different units.

A special case in which the Eckmann--Hilton argument is effective, is
any example in which $H_0\cong\K$. In this case,
$s=t=\epsilon_{H_1}\colon H_1\to\K$; $\imath=\eta_{H_1}\colon\K\to
H_1$; $H_1\pullback{s}{H_0}{t}H_1=H_1\otimes H_1$, and the
multiplication `$\circ$' is defined for all pairs of elements of
$H_1$. Furthermore, $\imath(1)\in H_1$ is a two-sided unit for
`$\circ$', \ie\,
\begin{equation}
  \circ(\imath(1)\otimes h)=h=\circ (h\otimes\imath(1)),
\end{equation}
for all $h\in H_1$. The Eckmann--Hilton argument then implies that for
all $h,h^\prime\in H_1$,
\begin{equation}
  \circ(h\otimes h^\prime) = \circ(h^\prime\otimes h) = hh^\prime =
  h^\prime h.
\end{equation}
Obviously, $H_1$ needs to be not only cocommutative, but also
commutative in order to admit such an example.
\end{remark}

\begin{remark}
By applying Proposition~\ref{prop_vertical} to
$\mathcal{C}=\cocCo_\K$, we see that in each strict cocommutative
trialgebra $(H_0,H_1,s,t,\imath,\circ)$, the bialgebra homomorphism
$\circ\colon H_1\pullback{s}{H_0}{t}H_1\to H_1$ is uniquely determined
by the remaining data. In Sweedler's notation, writing $\Delta(h)=\sum
h^{(1)}\otimes h^{(2)}\in H_1\otimes H_1$ for $h\in H_1$, it reads,
\begin{equation}
  \circ(h\otimes\tilde h)=\sum h^{(1)}S(\imath(s(h^{(2)})))\tilde h,
\end{equation}
for all $h\otimes\tilde h\in H_1\pullback{s}{H_0}{t}H_1$. The analogue
of vertical inversion is the bialgebra homomorphism $\US\colon H_1\to
H_1, h\mapsto\sum\imath(s(h^{(1)}))S(h^{(2)})\imath(t(h^{(3)}))$ which
satisfies for all $h\in H_1$,
\begin{eqnarray}
  t(\US(h))&=&s(h),\\
  s(\US(h))&=&t(h),\\
  \sum\circ(\US(h^{(1)})\otimes h^{(2)})&=&\imath(s(h)),\\
  \sum\circ(h^{(1)}\otimes\US(h^{(2)}))&=&\imath(t(h)),
\end{eqnarray}
and for all $h\otimes\tilde h\in H_1\pullback{s}{H_0}{t}H_1$,
\begin{equation}
  \US(\circ(h\otimes\tilde h))=\circ(\US(\tilde h)\otimes\US(h)).
\end{equation}
\end{remark}

%
\section{Commutative cotrialgebras}
%
\label{sect_comcotrialg}

In this section, we develop the concept dual to cocommutative
trialgebras. This construction is based on the functor
$\K(-)\colon\fGrp\to\comHopf_\K^\op$ which assigns to each finite
group $G$ its commutative Hopf algebra
$(\K(G),\mu,\eta,\Delta,\epsilon,S)$ of functions $\K(G)=\{\,f\colon
G\to\K\,\}$. This forms an associative unital algebra
$(\K(G),\mu,\eta)$ under pointwise operations in $\K$ and a Hopf
algebra using the operations inherited from the group structure of
$G$. The functor $\K(-)$ associates with each group homomorphism
$\phi\colon G\to H$ the bialgebra homomorphism
$\K(\phi)\colon\K(H)\to\K(G), f\mapsto \phi f$ (first $\phi$, then
$f$).

By internalization, we obtain a $2$-functor $\Cat(\K(-))$ which
associates with each finite $2$-group $G$ an internal cocategory in
the category of commutative Hopf algebras. We call such a structure a
\emph{commutative cotrialgebra}.  

\subsection{Finite colimits in the category of commutative Hopf algebras}
\label{sect_colimits}

We first recall the construction of finite colimits in the category
$\comHopf_\K$. The proofs of the following propositions are again
elementary.

\begin{proposition}
The initial object in any of the the categories $\Hopf_\K$,
$\comHopf_\K$, $\comBi_\K$ and $\comAlg_\K$ is given by $\K$
itself. For each object $A$, the unique morphism $\K\to A$ is the unit
operation $\eta\colon\K\to A$ of the underlying algebra.
\end{proposition}

\begin{proposition}
\label{prop_productcomhopf}
The binary coproduct in the category $\comHopf_\K$ is the tensor
product of Hopf algebras, \ie\ for commutative Hopf algebras $A$ and
$B$,
\begin{equation}
  A\textstyle\coprod B=A\otimes B,
\end{equation}
with the Hopf algebra homomorphisms $\imath_1=\id_A\otimes\eta_B\colon
A\to A\otimes B$ and $\imath_2=\eta_A\otimes\id_B\colon B\to A\otimes
B$. Here we have again suppressed the isomorphisms $A\otimes\K\cong A$
and $\K\otimes B\cong B$.

For each commutative Hopf algebra $D$ with homomorphisms $f_1\colon
A\to D$ and $f_2\colon B\to D$, there is a unique Hopf algebra
homomorphism $\phi\colon A\otimes B\to D$ such that $\imath_1\phi=f_1$
and $\imath_2\phi=f_2$. It is given by $\phi:=(f_1\otimes f_2)\mu$.
\end{proposition}

\begin{remark}
The binary coproduct in the categories $\comBi_\K$ and $\comAlg_\K$ is
given precisely by the same construction: the tensor product of
bialgebras or coalgebras, respectively. In the category $\Hopf_\K$ of
all (not necessarily commutative) Hopf algebras, the coproduct is in
general \emph{not} the tensor product of Hopf algebras. Commutativity
is essential in Proposition~\ref{prop_productcomhopf} in order to show
that $\phi\colon A\otimes B\to D$ is indeed a homomorphism of
algebras.
\end{remark}

\begin{proposition}
The binary coequalizer in the category $\Hopf_\K$ is given by a
quotient algebra. For Hopf algebras $C$, $D$ with homomorphisms
$f_1,f_2\colon C\to D$, the coequalizer of $f_1$ and $f_2$ is the
quotient algebra $D/I$ where $I$ is the two-sided (algebra) ideal
generated by $(f_1-f_2)(C)$, \ie\ by all elements of the form
$f_1(c)-f_2(c)$, $c\in C$. The associated colimiting cone is the
canonical projection $\pi\colon D\to D/I$. $I$ turns out to be a Hopf
ideal, and $\pi$ a homomorphism of Hopf algebras.

For each Hopf algebra $H$ with a homomorphism $\psi\colon D\to H$ such
that $f_1\psi=f_2\psi$, $\psi$ vanishes on $I\subseteq D$ so that it
descends to the quotient $D/I$ and gives rise to a Hopf algebra
homomorphism $\phi\colon D/I\to H$. This $\phi$ is the unique
Hopf algebra homomorphism for which $\psi=\pi\phi$.
\end{proposition}

\begin{remark}
The coequalizer in $\comHopf_\K$, $\comBi_\K$ and $\comAlg_\K$ is
given precisely by the same construction.
\end{remark}

\begin{corollary}
\label{corr_pushoutcommutativehopf}
The pushout in the category $\comHopf_\K$ is constructed as
follows. Let $A$, $B$, $C$ be commutative Hopf algebras with
homomorphisms $\sigma\colon C\to A$ and $\tau\colon C\to B$. The
pushout $A\pushout{\sigma}{C}{\tau}B$ is the quotient $(A\otimes B)/I$
where $I$ is the two-sided (algebra) ideal generated by
$\Phi(C)$. Here $\Phi$ denotes the linear map
$\Phi=(\sigma\otimes\eta_B-\eta_A\otimes\tau)\colon C\to A\otimes
B$. The colimiting cone is given by the Hopf algebra homomorphisms
$\imath_1=\id_A\otimes\eta_B\colon A\to A\pushout{\sigma}{C}{\tau} B$
and $\imath_2=\eta_A\otimes\id_A\colon B\to A\pushout{\sigma}{C}{\tau}
B$.

If $H$ is a commutative Hopf algebra with homomorphisms $f_1\colon
A\to H$ and $f_2\colon B\to H$ such that $\sigma f_1=\tau f_2$, then
there is a unique homomorphism $\phi\colon
A\pushout{\sigma}{C}{\tau}B\to H$ such that $\imath_1\phi=f_1$ and
$\imath_2\phi=f_2$. It is given by taking the quotient of the domain
of the homomorphism $\phi:=(f_1\otimes f_2)\mu\colon A\otimes B\to H$
modulo the ideal $I\subseteq A\otimes B$.
\end{corollary}

\subsection{Finite-limit preservation of the function algebra functor}

\begin{theorem}
The function algebra functor $\K(-)\colon\fGrp\to\comHopf_\K^\op$
preserves finite limits.
\end{theorem}

\begin{proof}
We verify in a direct calculation that $\K(-)$ preserves the terminal
object, binary products and binary equalizers. Since $\K(-)$ maps
$\fGrp$ to the opposite category $\comHopf_\K^\op$, this means $\K(-)$
maps the terminal object of $\fGrp$ to the initial object of
$\comHopf_\K$, maps binary products of $\fGrp$ to binary coproducts of
$\comHopf_\K$ and similarly binary equalizers to binary coequalizers.
\begin{myenumerate}
\item
  For the terminal object, obviously $\K(\{e\})\cong\K$ form
  isomorphic Hopf algebras. 
\item
  For binary products, choose for each finite group $G$, the basis
  ${\{\delta_g\}}_{g\in G}$ of $\K(G)$ where $\delta_g(g^\prime)=1$ if
  $g=g^\prime$ and $\delta_g(g^\prime)=0$ otherwise. For finite groups
  $G$ and $H$, there is an isomorphism of Hopf algebras,
\begin{equation}
  \K(G\times H)\to\K(G)\otimes\K(H),\quad
  \delta_{(g,h)}\to\delta_g\otimes\delta_h.
\end{equation}
  Let $p_1\colon G\times H\to G,(g,h)\mapsto g$ be the first
  projection. Then $\K(p_1)\colon\K(G)\to\K(G)\otimes\K(H)$ satisfies
  $\K(p_1)[\delta_g]=p_1\delta_g=((g^\prime,h^\prime)\mapsto
  \delta_g(p_1(g^\prime,h^\prime))=\delta_g(g^\prime))=(\id_{\K(G)}\otimes\eta_{\K(H)})(\delta_g)$
  for all $g\in G$. The analogous result holds for $p_2$. Therefore,
  $\K(-)$ maps the binary product of $\fGrp$ to the binary coproduct
  of $\comHopf_\K$ with the associated cones.
\item
  For equalizers, recall that the equalizer in $\fGrp$ of a pair of
  group homomorphisms $f_1,f_2\colon G\to H$ is the
  subgroup~\eqref{eq_equalizergrp} with its inclusion map. There is an
  isomorphism of Hopf algebras,
\begin{equation}
  \K(\eq(f_1,f_2))\cong\K(G)/I,
\end{equation}
where $I$ is the two-sided (algebra) ideal generated by all elements
of the form $\K(f_1)[f]-\K(f_2)[f]$ where $f\in\K(H)$. By choosing the
usual basis ${\{\delta_h\}}_{h\in H}$ of $\K(H)$, one sees that $I$ is
generated by all elements,
\begin{equation}
\label{eq_ideal1}
  \K(f_1)[\delta_h]-\K(f_2)[\delta_h]=f_1\delta_h-f_2\delta_h
  =\sum_{g\in G}\Bigl(\delta_h(f_1(g))-\delta_h(f_2(g))\Bigr)\delta_g,
\end{equation}
where $h\in H$. Exploiting the ideal property of $I\subseteq\K(G)$
and multiplying~\eqref{eq_ideal1} by all basis vectors
$\delta_{g^\prime}$, $g^\prime\in G$, of $\K(G)$, it follows that the
ideal $I$ is the following linear span,
\begin{equation}
  I=\Span_\K\Bigl\{\,\Bigl(\delta_h(f_1(g^\prime))-\delta_h(f_2(g^\prime))\Bigr)\delta_{g^\prime}
  \colon\quad h\in H,g^\prime\in G\,\Bigr\}.
\end{equation}
If therefore $f_1(g^\prime)\neq f_2(g^\prime)$, there is some
$h=f_1(g)$ such that the above generating set contains
$\delta_{g^\prime}$. If, however, $f_1(g^\prime)=f_2(g^\prime)$, then
$\delta_{g^\prime}$ is not contained. One can therefore read off a
basis for the quotient vector space $\K(G)/I$. It agrees with the
usual basis of $\K(\eq(f_1,f_2))$. It is furthermore obvious that
$\K(-)$ sends the inclusion map $\eq(f_1,f_2)\to G$ to the canonical
projection map $\K(G)\to \K(G)/I$ of $\K(G)$ onto the coequalizer of
$\K(f_1)$ and $\K(f_2)$.
\end{myenumerate}
\end{proof}

\begin{corollary}
The function algebra functor $\K(-)\colon\fGrp\to\comHopf_\K^\op$
preserves pullbacks. In particular, for groups $G$, $H$ and $K$ and
group homomorphisms $s\colon G\to K$ and $t\colon H\to K$, there is an
isomorphism of Hopf algebras,
\begin{equation}
  \K(G\times_K H)\cong\K(G)\textstyle\pushout{\K(s)}{\K(K)}{\K(t)}\K(H) 
    = (\K(G)\otimes\K(H))/I,
\end{equation}
where $I$ is the (algebra) ideal generated by all elements of the form
$\sigma(f)\otimes\eta_{\K(H)}-\eta_{\K(G)}\otimes\tau(f)$ for
$f\in\K(K)$ where we have written $\sigma=\K(s)$, $\tau=\K(t)$.
\end{corollary}

\begin{remark}
If the function algebra is viewed as a functor
$\K(-)\colon\fGrp\to\Hopf_k^\op$ into the category of all (not
necessarily commutative) Hopf algebras, it does not preserve
pullbacks. In fact, it does not even preserve binary products. This
can again be blamed on the inclusion functor $\comHopf_\K\to\Hopf_\K$
which does not preserve binary coproducts.
\end{remark}

\subsection{Definition of commutative cotrialgebras}
\label{sect_cotrialg}

\begin{definition}
The objects, morphisms and $2$-morphisms of the $2$-category,
\begin{equation}
  \comCoTri_\K:=\CoCat(\comHopf_\K),
\end{equation}
are called \emph{strict commutative cotrialgebras}, their
homomorphisms and $2$-homomorphisms, respectively.
\end{definition}

The main consequence of the preceding section is the existence of the
following $2$-functor which can be used in order to construct examples
of commutative cotrialgebras from strict finite $2$-groups and from
finite crossed modules.

\begin{proposition}
\label{prop_cotrialg}
The functor $\K(-)\colon\fGrp\to\comHopf_\K$ gives rise to a
$2$-functor,
\begin{equation}
  \Cat(\K(-))\colon\fTGrp\to\comCoTri_\K.
\end{equation}
\end{proposition}

Let us finally unfold the definition of a strict commutative
cotrialgebra. 

\begin{proposition}
A strict commutative cotrialgebra
$H=(H_0,H_1,\sigma,\tau,\Uepsilon,\UDelta)$ consists of commutative
Hopf algebras $H_0$ and $H_1$ over $\K$ and of bialgebra homomorphisms
$\sigma,\tau\colon H_0\to H_1$, $\Uepsilon\colon H_1\to H_0$ and
$\UDelta\colon H_1\to H_1\pushout{\sigma}{H_0}{\tau}H_1$ such that the
diagrams dual to~\eqref{eq_intcat1}--\eqref{eq_intcat4} hold, renaming
$C_j\mapsto H_h$; $s\mapsto\sigma$, $t\mapsto\tau$;
$\imath\mapsto\Uepsilon$ and $\circ\mapsto\UDelta$.
\end{proposition}

\begin{remark}
The terminology \emph{cotrialgebra} originates from the three
operations that are defined on the vector space $H_1$. Here we have
again the multiplication and the comultiplication of the Hopf algebra
$H_1$, but with an additional, partially defined comultiplication,
\begin{equation}
  \UDelta\colon H_1\to H_1\textstyle\pushout{\sigma}{H_0}{\tau}H_1.
\end{equation}
The notion of a partially defined comultiplication is precisely dual
to that of the partially defined multiplication in a trialgebra
(Remark~\ref{rem_trialgebra}). The Hopf algebra homomorphism does not
map into $H_1\otimes H_1$, but rather into a suitable quotient of this
algebra which is given by the pushout
(Corollary~\ref{corr_pushoutcommutativehopf}). Similarly, the
partially defined comultiplication has a local counit,
\begin{equation}
  \Uepsilon\colon H_1\to H_0.
\end{equation}
In particular, the diagram dual to~\eqref{eq_intcat3},
\begin{equation}
\begin{aligned}
\xymatrix{
  H_1\ar[rr]^{\UDelta}\ar[dd]_{\UDelta}&&
    H_1\pushout{\sigma}{H_0}{\tau}H_1\ar[dd]^{{[\id_{H_1};\UDelta]}_{\id_{H_0}}}\\
  \\
  H_1\pushout{\sigma}{H_0}{\tau}H_1\ar[rr]_{{[\UDelta;\id_{H_1}]}_{\id_{H_0}}}&&
  H_1\pushout{\sigma}{H_0}{\tau}H_1\pushout{\sigma}{H_0}{\tau}H_1
}
\end{aligned}
\end{equation}
states the coassociativity for the partially defined comultiplication
while the diagram dual to~\eqref{eq_intcat4},
\begin{equation}
\begin{aligned}
\xymatrix{
  &&H_1\ar[lldd]_{\imath_1}\ar[dd]^{\UDelta}\ar[rrdd]^{\imath_2}\\
  \\
  H_1\pushout{\sigma}{H_0}{\id_{H_0}}H_0&&
    H_1\pushout{\sigma}{H_0}{\tau}H_1\ar[ll]^{{[\id_{H_1};\Uepsilon]}_{\id_{H_0}}}
      \ar[rr]_{{[\Uepsilon;\id_{H_1}]}_{\id_{H_0}}}&&
    H_0\pushout{\id_{H_0}}{H_0}{\tau}H_1
}
\end{aligned}
\end{equation}
states its counit property. As far as the compatibility of the three
operations on $H_1$ and the avoidance of the Eckmann--Hilton argument
is concerned, everything is dual to the case of a cocommutative
trialgebra.
\end{remark}

\subsection{Compact topological $2$-groups and commutative $C^\ast$-cotrialgebras}
\label{sect_top2grp}

Whereas the definition of the group algebra $\K[G]$
(Section~\ref{sect_coctrialg}) is available for any group $G$, we have
to refine the definition of the function algebra $\K(G)$
(Section~\ref{sect_comcotrialg}) if $G$ is supposed to be more general
than a finite group\footnote{This is because we want to provide
$\K(G)$ with a comultiplication
$\Delta\colon\K(G)\to\K(G)\otimes\K(G)$ which is induced by the group
multiplication $G\times G\to G$, \ie\ we need an isomorphism
$\K(G\times G)\cong\K(G)\otimes\K(G)$. The algebraic tensor product,
however, may be insufficient for this.}.

One possibility is to restrict the function algebra to the algebraic,
\ie\ polynomial functions. If $G$ is a compact topological group,
however, a good choice of function algebra is the $C^\ast$-algebra
$C(G)$ of continuous complex-valued functions on $G$. The resulting
theory is particularly powerful because the $C^\ast$-algebra $C(G)$
contains the full information in order to reconstruct $G$ using
Gel'fand representation theory. This fact can be expressed as an
equivalence between the category of compact topological groups and a
suitable category of $C^\ast$-algebras. In Appendix~\ref{app_topgrp},
we review the construction of these categories and how to derive their
equivalence from Gel'fand representation theory. In the following, we
just summarize the relevant results and generalize the preceding
subsection, in particular Proposition~\ref{prop_cotrialg}, to the case
of compact topological groups.

\begin{definition}
\begin{myenumerate}
\item
  A \emph{compact topological group} $(G,\mu,\eta,\zeta)$ is a compact
  Hausdorff space $G$ with continuous maps,
\begin{alignat}{2}
\mu&\colon   G\times G\to G,&\qquad&\mbox{(multiplication)},\\
\eta&\colon  \{e\}\to G,    &\qquad&\mbox{(unit)},\\
\zeta&\colon G\to G,        &\qquad&\mbox{(inversion)},
\end{alignat}
  such that $(\mu\times\id_G)\mu=(\id_G\times\mu)\mu$,
  $(\eta\times\id_G)\mu=p_2$, $(\id_G\times\eta)\mu=p_1$,
  $\delta(\id_G\times\zeta)\mu=t\eta$ and
  $\delta(\zeta\times\id_G)\mu=t\eta$. Here $p_1\colon G\times G\to G,
  (g_1,g_2)\mapsto g_1$ is the projection onto the first factor,
  similarly $p_2$ onto the second factor; $\delta\colon G\to G\times
  G$, $g\mapsto (g,g)$ is the diagonal map, and $t\colon G\to\{e\}$,
  $g\to e$.
\item
  Let $(G,\mu,\eta,\zeta)$ and
  $(G^\prime,\mu^\prime,\eta^\prime,\zeta^\prime)$ be compact
  topological groups. A \emph{homomorphism of compact topological
  groups} $f\colon G\to G^\prime$ is a continuous map for which $\mu
  f=(f\times f)\mu^\prime$, $\eta f=\eta^\prime$ and $\zeta
  f=f\zeta^\prime$.
\item
  There is a category $\compTopGrp$ whose objects are compact
  topological groups and whose morphisms are homomorphisms of compact
  topological groups.
\end{myenumerate}
\end{definition}

\begin{definition}
\label{def_hopfcstar}
\begin{myenumerate}
\item
  A \emph{commutative Hopf $C^\ast$-algebra} $(H,\Delta,\epsilon,S)$
  is a commutative unital $C^\ast$-algebra $H$ with unital
  $\ast$-homomorphisms,
\begin{alignat}{2}
\Delta&\colon   H\to H\otimes_\ast H,&\qquad&\mbox{(comultiplication)},\\
\epsilon&\colon H\to \C,             &\qquad&\mbox{(counit)},\\
S&\colon        H\to H,              &\qquad&\mbox{(antipode)},
\end{alignat}
  such that $\Delta(\Delta\otimes\id_H)=\Delta(\id_H\otimes\Delta)$,
  $\Delta(\id_H\otimes\epsilon)=\id_H\otimes 1_H$,
  $\Delta(\epsilon\otimes\id_H)=1_H\otimes\id_H$,
  $\Delta(S\otimes\id_H)\mu=\epsilon\eta$ and $\Delta(\id_H\otimes
  S)\mu=\epsilon\eta$ (composition is read from left to right). Here
  $1_H$ is the unit of $H$, $\eta\colon\C\to H$, $1\mapsto 1_H$, and
  $\mu\colon H\otimes_\ast H\to H$ denotes multiplication in $H$. The
  notation $\otimes_\ast$ refers to the unique completion of the
  algebraic tensor product to a $C^\ast$-algebra (see
  Appendix~\ref{app_topgrp}).
\item
  Let $(H,\Delta,\epsilon,S)$ and
  $(H^\prime,\Delta^\prime,\epsilon^\prime,S^\prime)$ be commutative
  Hopf $C^\ast$-algebras. A \emph{homomorphism of commutative Hopf
  $C^\ast$-algebras} is a unital $\ast$-homomorphism $f\colon H\to
  H^\prime$ for which $\Delta(f\otimes f)=f\Delta^\prime$,
  $\epsilon=f\epsilon^\prime$ and $Sf=fS^\prime$.
\item
  There is a category $\comHopfCstAlg$ whose objects are commutative
  Hopf $C^\ast$-algebras and whose morphisms are homomorphisms of
  commutative Hopf $C^\ast$-algebras.
\end{myenumerate}
\end{definition}

\begin{remark}
\begin{myenumerate}
\item
  In a generic Hopf algebra, one usually defines the antipode $S\colon
  H\to H$ as a linear map with the property that
  $\Delta(S\otimes\id_H)\mu=\epsilon\eta$ and $\Delta(\id_H\otimes
  S)\mu=\epsilon\eta$. It then turns out that $S$ is an algebra
  \emph{anti}-homomorphism. Definition~\ref{def_hopfcstar} looks a bit
  strange because it defines $S$ as an algebra homomorphism (which is
  here, of course, the same as an anti-homomorphism since $H$ is
  commutative). The reason for this choice is our derivation of
  Theorem~\ref{thm_usegelfand} in Appendix~\ref{app_topgrp}.
\item
  Note that the terminology \emph{comultiplication} is used by some
  authors only if $\Delta$ maps into the algebraic tensor product
  $H\otimes H$. Our comultiplication which maps into a completion of
  the tensor product would then be called a \emph{topological
  comultiplication}.
\end{myenumerate}
\end{remark}

There is a functor $C(-)\colon\compTopGrp\to\comHopfCstAlg^\op$ which
assigns to each compact topological group $G$ the commutative unital
$C^\ast$-algebra $C(G)$ of continuous complex-valued functions on
$G$. Using the group operations of $G$, this $C^\ast$-algebra is
equipped with the structure of a Hopf $C^\ast$-algebra by setting
$(\Delta f)[g_1,g_2]:=f(\mu(g_1,g_2))$, $\epsilon f:=f(\eta(1))$ and
$(Sf)[g]:=f(\zeta(g))$. There is also a functor
$\sigma(-)\colon\comHopfCstAlg^\op\to\compTopGrp$ which assigns to
each commutative Hopf $C^\ast$-algebra $H$ its Gel'fand spectrum
$\sigma(H)$ which is a compact topological group using the coalgebra
structure of $H$. For more details, we refer to
Appendix~\ref{app_topgrp} in which we also review how to derive the
following theorem.

\begin{theorem}
\label{thm_usegelfand}
There is an equivalence of categories
$\compTopGrp\simeq\comHopfCstAlg^\op$ provided by the functors
\begin{equation}
  C(-)\colon\compTopGrp\to\comHopfCstAlg^\op
\end{equation}
and
\begin{equation}
  \sigma(-)\colon\comHopfCstAlg^\op\to\compTopGrp.
\end{equation}
\end{theorem}

The category $\compTopGrp$ is finitely complete. In particular, the
terminal object is the trivial group $\{e\}$, the binary product
$G\prod H=G\times H$ of two compact topological groups $G$ and $H$ is
the direct product with the product topology and the usual
projections, and the equalizer of a pair of homomorphisms of compact
topological groups $f_1,f_2\colon G\to H$ is the (closed) subgroup
$E=\{\,g\in G\colon f_1(g)=f_2(g)\,\}\subseteq G$ with the induced
topology and the canonical inclusion map.

Since the functors $C(-)$ and $\sigma(-)$ form an equivalence of
categories, they both preserve finite limits and finite colimits, and
$\comHopfCstAlg$ is finitely cocomplete.

\begin{definition}
The objects of the $2$-category,
\begin{equation}
  \compTopTGrp := \Cat(\compTopGrp)
\end{equation}
are called \emph{strict compact topological $2$-groups} and the
objects of the $2$-category,
\begin{equation}
  \comCstCoTri := \CoCat(\comHopfCstAlg)
\end{equation}
\emph{strict commutative $C^\ast$-cotrialgebras}.
\end{definition}

We can immediately employ the general theory and invoke
Corollary~\ref{corr_equivalence} in order to obtain the following
result which generalizes Proposition~\ref{prop_cotrialg} without the
need to verify limit preservation of the functors by hand.

\begin{theorem}
\label{thm_equivcst}
There is a $2$-equivalence between the $2$-categories
\begin{equation}
\label{eq_equivcst}
  \compTopTGrp\simeq\comCstCoTri
\end{equation}
provided by the functors
\begin{equation}
  \Cat(C(-))\colon\compTopTGrp\to\comCstCoTri
\end{equation}
and 
\begin{equation}
  \Cat(\sigma(-))\colon\comCstCoTri\to\compTopTGrp.
\end{equation}
\end{theorem}

We can use the $2$-functor $\Cat(C(-))$ in order to construct a strict
commutative $C^\ast$-cotrialgebra $H:=\Cat(C(-))[G]$ from each strict
compact topological $2$-group $G$. Conversely, we can recover the
strict compact topological $2$-group $G$ from the strict commutative
$C^\ast$-cotrialgebra $H$ up to isomorphism as
$G\cong\Cat(\sigma(-))[H]$. The $2$-functor $\Cat(\sigma(-))$ is
therefore the appropriate generalization of the Gel'fand spectrum to
the case of commutative $C^\ast$-cotrialgebras and
Theorem~\ref{thm_equivcst} the desired generalization of Gel'fand
representation theory.

The general remarks of the Sections~\ref{sect_colimits}
to~\ref{sect_cotrialg} generalize from finite groups to compact
topological groups if we replace the algebra $\K(G)$ by $C(G)$.

\begin{remark}
\label{rem_vertical2}
Applying Proposition~\ref{prop_vertical} to any of the categories
$\mathcal{C}=\comAlg_\K^\op$ or $\mathcal{C}=\comUnCstAlg^\op$, we see
that in any strict commutative [$C^\ast$-]cotrialgebra
$(H_0,H_1,\sigma,\tau,\Uepsilon,\UDelta)$, the homomorphism of Hopf
[$C^\ast$-]algebras $\UDelta\colon H_1\to
H_1\pushout{\sigma}{H_0}{\tau}H_1$ is uniquely determined by,
\begin{equation}
  \UDelta(h) = \sum h^{(1)}\sigma(\Uepsilon(S(h^{(2)})))\otimes h^{(3)},
\end{equation}
for all $h\in H_1$, where we denote by $\Delta(h)=\sum h^{(1)}\otimes
h^{(2)}$ the comultiplication of $H_1$. There is a homomorphism of
Hopf [$C^\ast$-]algebras,
\begin{equation}
\label{eq_antipode2}
  \US\colon H_1\to H_1,\quad 
  h\mapsto\sum\sigma(\Uepsilon(h^{(1)}))S(h^{(2)})\tau(\Uepsilon(h^{(3)})),
\end{equation}
such that,
\begin{eqnarray}
  \tau\US&=&\sigma,\\
  \sigma\US&=&\tau,\\
  \UDelta(\US\otimes\id_{H_1})\mu_{H_1}&=&\Uepsilon\sigma,\\
  \UDelta(\id_{H_1}\otimes\US)\mu_{H_1}&=&\Uepsilon\tau,\\
  \UDelta(\US\otimes\US)&=&\US\UDelta^\op,
\end{eqnarray}
where composition is read from left to right.
\end{remark}

%
\section{Symmetric Hopf categories}
%
\label{sect_hopfcat}

\subsection{Tannaka--Kre\v\i n duality for commutative Hopf algebras}

In order to define and study symmetric Hopf categories, we proceed in
close analogy to our treatment of Gel'fand representation theory in
Section~\ref{sect_top2grp}. We start from the notion of a strict
commutative cotrialgebra (Section~\ref{sect_cotrialg}) which is
defined as an internal cocategory in the category of commutative Hopf
algebras over some field $\K$. We then exploit the equivalence of the
category of commutative Hopf algebras with a suitably chosen category
whose objects are symmetric monoidal categories. This result is known
as Tannaka--Kre\v\i n duality and allows us to reconstruct the
commutative Hopf algebra from its category of finite-dimensional
comodules. Both functors involved in the equivalence preserve
colimits. We can therefore make use of the $2$-functor $\Cat(-)$ in
order to define symmetric Hopf categories and to establish a
generalization of Tannaka--Kre\v\i n duality between commutative
cotrialgebras and symmetric Hopf categories.

Standard references on Tannaka--Kre\v\i n duality
are~\cite{Sa72,DeMi82,JoSt91}. We are aiming for an equivalence of
categories, and it is difficult to find such a result explicitly
stated in the literature. We follow the presentation of~\cite{Sc92}
which comes closest to our goal. In Appendix~\ref{app_tkduality}, we
summarize how the results of~\cite{Sc92} can be employed in order to
establish the desired equivalence of categories. In the present
section, we just state the relevant definitions and results. We use
the term `symmetric monoidal category' for a tensor category and
`braided monoidal category' for a quasi-tensor category
following~\cite{Ma95}.

\begin{definition}
\label{def_catvect}
Let $\mathcal{V}$ be a category.
\begin{myenumerate}
\item
  A \emph{category over} $\mathcal{V}$ is a pair
  $(\mathcal{C},\mathcal{\omega})$ of a category $\mathcal{C}$ and a
  functor $\omega\colon\mathcal{C}\to\mathcal{V}$.
\item
  Let $(\mathcal{C},\omega)$ and $(\mathcal{C}^\prime,\omega^\prime)$
  be categories over $\mathcal{V}$. A \emph{functor over}
  $\mathcal{V}$ is an equivalence class $[F,\zeta]$ of pairs
  $(F,\zeta)$ where $F\colon\mathcal{C}\to\mathcal{C}^\prime$ is a
  functor and $\zeta\colon\omega\Rightarrow F\cdot\omega^\prime$ is a
  natural isomorphism. Two such pairs $(F,\zeta)$ and $(\tilde
  F,\tilde\zeta)$ are considered equivalent if and only if there is a
  natural isomorphism $\phi\colon F\Rightarrow\tilde F$ such that,
\begin{equation}
  \tilde\zeta = (\phi\cdot\id_{\omega^\prime})\circ\zeta.
\end{equation}
  Recall that vertical `$\circ$' composition of natural
  transformations is read from right to left, but horizontal `$\cdot$'
  composition from left to right.
\item
  For each category $(\mathcal{C},\omega)$ over $\mathcal{V}$, we
  define the \emph{identity functor}
  $\id_{(\mathcal{C},\omega)}:=[1_\mathcal{C},\id_\omega]$ over
  $\mathcal{V}$.
\item
  Let $(\mathcal{C},\omega)$, $(\mathcal{C}^\prime,\omega^\prime)$ and
  $(\mathcal{C}^{\prime\prime},\omega^{\prime\prime})$ be categories
  over $\mathcal{V}$ and
  $[F,\zeta]\colon(\mathcal{C},\omega)\to(\mathcal{C}^\prime,\omega^\prime)$
  and
  $[G,\theta]\colon(\mathcal{C}^\prime,\omega^\prime)\to(\mathcal{C}^{\prime\prime},\omega^{\prime\prime})$
  be functors over $\mathcal{V}$. Then their \emph{composition} is
  defined as,
\begin{equation}
  [F,\zeta]\cdot[G,\theta] := [F\cdot G,(\id_F\cdot\theta)\circ\zeta].
\end{equation}
\end{myenumerate}
\end{definition}

\begin{definition}
Let $\mathcal{V}$ be a monoidal category.
\begin{myenumerate}
\item
  A \emph{monoidal category over} $\mathcal{V}$ is a category
  $(\mathcal{C},\omega)$ over $\mathcal{V}$ such that $\mathcal{C}$ is
  a monoidal category and $\omega\colon\mathcal{C}\to\mathcal{V}$ a
  monoidal functor.
\item
  Let $(\mathcal{C},\omega)$ and $(\mathcal{C}^\prime,\omega^\prime)$
  be monoidal categories over $\mathcal{V}$. A \emph{monoidal functor
  over} $\mathcal{V}$ is a functor
  $[F,\zeta]\colon(\mathcal{C},\omega)\to(\mathcal{C}^\prime,\omega^\prime)$
  over $\mathcal{V}$ such that
  $F\colon\mathcal{C}\to\mathcal{C}^\prime$ is a monoidal functor and
  $\zeta$ is a monoidal natural isomorphism.
\end{myenumerate}
\end{definition}

\begin{definition}
Let $\mathcal{V}$ be a rigid monoidal category.
\begin{myenumerate}
\item
  A \emph{rigid monoidal category over} $\mathcal{V}$ is a monoidal
  category $(\mathcal{C},\omega)$ over $\mathcal{V}$ such that
  $\mathcal{C}$ is a rigid monoidal category (\ie\ every object has
  got a left-dual).
\item
  Let $(\mathcal{C},\omega)$ and $(\mathcal{C}^\prime,\omega^\prime)$
  be rigid monoidal categories over $\mathcal{V}$. A \emph{rigid
  monoidal functor over} $\mathcal{V}$ is a monoidal functor
  $[F,\zeta]\colon(\mathcal{C},\omega)\to(\mathcal{C}^\prime,\omega^\prime)$
  over $\mathcal{V}$ such that the functor
  $F\colon\mathcal{C}\to\mathcal{C}^\prime$ preserves left-duals.
\end{myenumerate}
\end{definition}

\begin{definition}
Let $\mathcal{V}$ be a symmetric monoidal category.
\begin{myenumerate}
\item
  A \emph{symmetric monoidal category over} $\mathcal{V}$ is a
  monoidal category $(\mathcal{C},\omega)$ over $\mathcal{V}$ such
  that $\mathcal{C}$ is a symmetric monoidal category and $\omega$ is
  a symmetric monoidal functor.
\item
  Let $(\mathcal{C},\omega)$ and $(\mathcal{C}^\prime,\omega^\prime)$
  be symmetric monoidal categories over $\mathcal{V}$. A
  \emph{symmetric monoidal functor over} $\mathcal{V}$ is a monoidal
  functor
  $[F,\zeta]\colon(\mathcal{C},\omega)\to(\mathcal{C}^\prime,\omega^\prime)$
  over $\mathcal{V}$ such that for all objects $X,Y$ of $\mathcal{C}$
  the following diagram commutes,
\begin{equation}
\begin{aligned}
\label{eq_symmetricfunctor}
\xymatrix{
  \omega^\prime(F(X)\otimes F(Y))\ar[rr]^{\omega^\prime(\xi_{X,Y})}
    \ar[dd]_{\omega^\prime(\sigma^\prime_{F(X),F(Y)})}
    &&\omega^\prime(F(X\otimes Y))\ar[dd]^{\omega^\prime(F(\sigma_{X,Y}))}\\
  \\
  \omega^\prime(F(Y)\otimes F(X))\ar[rr]_{\omega^\prime(\xi_{Y,X})}
    &&\omega^\prime(F(Y\otimes X)).
}
\end{aligned}
\end{equation}
  Here $\sigma_{-,-}$ denotes the (symmetric) braiding of
  $\mathcal{C}$, $\sigma^\prime_{-,-}$ the (symmetric) braiding of
  $\mathcal{C}^\prime$, and $\xi_{-,-}$ the natural equivalence that
  determines the monoidal structure of the functor
  $F\colon\mathcal{C}\to\mathcal{C}^\prime$. The
  diagram~\eqref{eq_symmetricfunctor} is the image under the functor
  $\omega^\prime\colon\mathcal{C}^\prime\to\mathcal{V}$ of the
  condition that $F$ be a symmetric monoidal functor.
\end{myenumerate}
\end{definition}

\begin{proposition}
\label{prop_rigid}
Let $\mathcal{V}$ be a rigid monoidal category, $(\mathcal{C},\omega)$
and $(\mathcal{C}^\prime,\omega^\prime)$ be rigid monoidal categories
over $\mathcal{V}$ and
\begin{equation}
  [F,\zeta]\colon(\mathcal{C},\omega)\to(\mathcal{C}^\prime,\omega^\prime)
\end{equation}
be a monoidal functor over $\mathcal{V}$. Then $[F,\zeta]$ is a rigid
monoidal functor over $\mathcal{V}$.
\end{proposition}

\begin{proof}
Standard, see, for example~\cite{DeMi82}.
\end{proof}

\begin{proposition}
\label{prop_symmetric}
Let $\mathcal{V}$ be a symmetric monoidal category,
$(\mathcal{C},\omega)$ and $(\mathcal{C}^\prime,\omega^\prime)$ be
symmetric monoidal categories over $\mathcal{V}$ and
$[F,\zeta]\colon(\mathcal{C},\omega)\to(\mathcal{C}^\prime,\omega^\prime)$
be a monoidal functor over $\mathcal{V}$. Then $[F,\zeta]$ is a
symmetric monoidal functor over $\mathcal{V}$.
\end{proposition}

\begin{proof}
The proof involves a huge commutative diagram, but is otherwise an
immediate consequence of the definitions.
\end{proof}

\begin{definition}
Let $\K$ be a field. A \emph{reconstructible category over} $\Vect_\K$
is a rigid symmetric monoidal category $(\mathcal{C},\omega)$ over
$\Vect_\K$ such that $\mathcal{C}$ is a $\K$-linear Abelian
essentially small category and $\omega\colon\mathcal{C}\to\Vect_\K$ is
an exact faithful $\K$-linear functor into $\fdVect_\K$. Recall that
$k$-linearity requires the tensor product of $\mathcal{C}$ to be
$k$-bilinear and the braiding to be $k$-linear. A
\emph{reconstructible functor over} $\Vect_\K$ is any monoidal functor
over $\Vect_\K$. $\mathfrak{C}^{\otimes,\ast,s}_{\K,rec}$ denotes the
category whose objects are reconstructible categories over $\Vect_\K$
and whose morphisms are reconstructible functors over $\Vect_\K$.
\end{definition}

The following theorem states the Tannaka--Kre\v\i n duality between
commutative Hopf algebras over $\K$ and reconstructible categories
over $\Vect_\K$. For more details, see Appendix~\ref{app_tkduality}.

\begin{theorem}
\label{thm_statetkduality}
There is an equivalence of categories,
\begin{equation}
\label{eq_statetkduality}
  \mathfrak{C}^{\otimes,\ast,s}_{\K,rec}\simeq\comHopf_\K,
\end{equation}
given by the functors
\begin{equation}
  \coend(-)\colon\mathfrak{C}^{\otimes,\ast,s}_{\K,rec}\to\comHopf_\K
\end{equation}
and
\begin{equation}
  \comod(-)\colon\comHopf_\K\to\mathfrak{C}^{\otimes,\ast,s}_{\K,rec}
\end{equation}
The functor $\coend(-)$ associates with each reconstructible category
over $\Vect_\K$ the coendomorphism coalgebra over $\K$ which can be
shown to form a commutative Hopf algebra. The functor $\comod(-)$
associates with each commutative Hopf algebra $H$ over $\K$ the
category ($\mathcal{M}^H,\omega^H)$ over $\Vect_\K$ where
$\mathcal{M}^H$ is the category of finite-dimensional
right-$H$-comodules and $\omega^H\colon\mathcal{M}^H\to\Vect_\K$ is
the forgetful functor. This category over $\Vect_\K$ can be shown to
be reconstructible.
\end{theorem}

From the equivalence~\eqref{eq_statetkduality}, it follows that both
functors $\coend(-)$ and $\comod(-)$ preserve colimits. Since the
category $\comHopf_\K$ is finitely cocomplete
(Section~\ref{sect_colimits}), the same holds for the category
$\mathfrak{C}^{\otimes,\ast,s}_{\K,rec}$.

\subsection{Definition of symmetric Hopf categories}

\begin{definition}
\label{def_hopfcat}
The objects of the $2$-category
$\symHopf_\K:=\CoCat(\mathfrak{C}^{\otimes,\ast,s}_{\K,rec})$ are
called \emph{strict symmetric Hopf categories} over $\K$.
\end{definition}

This definition deviates from similar definitions used in the
literature~\cite{CrFr94,Ne99} in a subtle way. First, we have
restricted ourselves to the symmetric case. Otherwise, it would not be
possible to use the techniques of internalization. Second, in this
definition which we unfold in the subsequent proposition and remark,
we have a comultiplication functor which maps not into the external
tensor product, but rather into a pushout. Third, as we show in
Proposition~\ref{prop_antipode} below, a functorial antipode is always
present. We therefore use the term \emph{Hopf category} rather than
\emph{bialgebra category}, \emph{bimonoidal category} or
\emph{bitensor category}.

\begin{proposition}
A strict symmetric Hopf category
$\mathcal{C}=(\mathcal{C}_0,\mathcal{C}_1,\sigma,\tau,\Uepsilon,\UDelta)$
consists of reconstructible categories $\mathcal{C}_0$ and
$\mathcal{C}_1$ over $\Vect_\K$ and of reconstructible functors
$\sigma,\tau\colon\mathcal{C}_0\to\mathcal{C}_1$, $\Uepsilon\colon
\mathcal{C}_1\to\mathcal{C}_0$ and $\UDelta\colon\mathcal{C}_1\to
\mathcal{C}_1\pushout{\sigma}{\mathcal{C}_0}{\tau}\mathcal{C}_1$ such
that the diagrams dual to~\eqref{eq_intcat1}--\eqref{eq_intcat4} hold
(with the obvious renamings)
\end{proposition}

\begin{remark}
\label{rem_hopfcat}
The reason for calling this structure a \emph{strict symmetric Hopf
category} are the properties of the category $\mathcal{C}_1$. It is a
symmetric monoidal category over $\Vect_K$, and so $\mathcal{C}_1$ has
a symmetric tensor product which defines a functor,
\begin{equation}
\label{eq_hopftensor}
  \hat\otimes\colon\mathcal{C}_1\boxtimes\mathcal{C}_1\to\mathcal{C}_1,
\end{equation}
which forms the categorical analogue of a multiplication operation
(functorial antipode). Here $\mathcal{C}\boxtimes\mathcal{D}$ denotes
the external tensor product of two $\K$-linear Abelian categories
$\mathcal{C}$ and $\mathcal{D}$ over $\Vect_\K$. The $\K$-linear
Abelian category $\mathcal{C}\boxtimes\mathcal{D}$ over $\Vect_\K$
together with the $\K$-bilinear functor
$\boxtimes\colon\mathcal{C}\times\mathcal{D}\to\mathcal{C}\boxtimes\mathcal{D}$
over $\Vect_\K$ is determined by the universal property that every
$\K$-bilinear functor
$F\colon\mathcal{C}\times\mathcal{D}\to\mathcal{E}$ over $\Vect_\K$
factors though $\boxtimes$ with unique $\hat F$ (unique as a morphism
of $\mathfrak{C}^{\otimes,\ast,s}_{\K,rec}$), \ie\,
\begin{equation}
\begin{aligned}
\xymatrix{
  \mathcal{C}\times\mathcal{D}\ar[rr]^\boxtimes\ar[rrdd]_F&&
    \mathcal{C}\boxtimes\mathcal{D}\ar[dd]^{\hat F}\\ \\
  &&\mathcal{E}
}
\end{aligned}
\end{equation}
In particular, the tensor product
$\otimes\colon\mathcal{C}_1\times\mathcal{C}_1\to\mathcal{C}_1$ is a
$\K$-bilinear functor and therefore defines the
functor~\eqref{eq_hopftensor} by,
\begin{equation}
\begin{aligned}
\xymatrix{
  \mathcal{C}_1\times\mathcal{C}_1\ar[rr]^\boxtimes\ar[rrdd]_\otimes&&
    \mathcal{C}_1\boxtimes\mathcal{C}_1\ar[dd]^{\hat\otimes}\\ \\
  &&\mathcal{C}_1
}
\end{aligned}
\end{equation}
In addition to~\eqref{eq_hopftensor}, there is the functorial
comultiplication,
\begin{equation}
\label{eq_hopfcotensor}
  \UDelta\colon\mathcal{C}_1\to
    \mathcal{C}_1\pushout{\sigma}{\mathcal{C}_0}{\tau}\mathcal{C}_1
\end{equation}
The compatibility of $\UDelta$ with $\hat\otimes$ is encoded in the
requirement that $\UDelta$ be a reconstructible functor over
$\Vect_\K$, \ie\ in particular a monoidal functor over $\Vect_\K$. If
the functor $\UDelta$ mapped into the coproduct of the category
$\mathfrak{C}^{\otimes,\ast,s}_{\K,rec}$ (which is the external tensor
product) rather than into the pushout, the category $\mathcal{C}_1$
would have precisely the structure expected from a categorification of
the notion of a commutative bialgebra, namely~\eqref{eq_hopftensor}
and a functor over $\Vect_\K$,
\begin{equation}
  \UDelta\colon\mathcal{C}_1\to\mathcal{C}_1\boxtimes\mathcal{C}_1.
\end{equation}
In the terminology of Neuchl~\cite{Ne99}, such a category
$\mathcal{C}_1$ would be called a \emph{bimonoidal category} or
\emph{$2$-bialgebra}. In contrast to Neuchl's definition, however, our
functor $\UDelta$ of~\eqref{eq_hopfcotensor} does not map into the
coproduct, but rather into a pushout. This is in complete analogy to
the situation for $2$-groups and (co-)trialgebras where one of the
operations is always partially defined, which helps avoid the
Eckmann--Hilton argument.
\end{remark}

\subsection{Tannaka--Kre\v\i n reconstruction of commutative cotrialgebras}

Similarly to our treatment of commutative cotrialgebras, we can employ
the general theory and invoke Corollary~\ref{corr_equivalence} in
order to obtain the generalization of Tannaka--Kre\v\i n duality to
commutative cotrialgebras.

\begin{theorem}
There is a $2$-equivalence between the $2$-categories
\begin{equation}
\label{eq_twotkduality}
  \comCoTri_\K\simeq\symHopf_\K
\end{equation}
provided by the functors
\begin{equation}
  \Cat(\comod(-))\colon\comCoTri_\K\to\symHopf_\K
\end{equation}
and
\begin{equation}
  \Cat(\coend(-))\colon\symHopf_\K\to\comCoTri_\K.
\end{equation}
\end{theorem}

\begin{remark}
\label{rem_antipode}
We can therefore use the $2$-functor $\Cat(\comod(-))$ in order to
construct a strict symmetric Hopf category from each strict
commutative cotrialgebra. Conversely, the $2$-functor
$\Cat(\coend(-))$ allows us to reconstruct the cotrialgebra from the
Hopf category up to isomorphism.

Our definition of strict symmetric Hopf categories is engineered in
such a way that Tannaka--Kre\v\i n duality can be generalized and that
we obtain a large family of examples from strict $2$-groups. Comparing
our definition with the literature~\cite{CrFr94,Ne99}, there is one
major discrepancy. Our functorial
comultiplication~\eqref{eq_hopfcotensor} maps into a pushout rather
than into the external tensor product. This immediately raises the
question of how Neuchl's work~\cite{Ne99} can be extended to include
our notion of Hopf categories.

A further question concerns the categorical analogue of an
antipode. The work of Crane and Frenkel~\cite{CrFr94} and of Carter,
Kauffman and Saito~\cite{CaKa99} which aims for a categorification of
Kuperberg's invariant~\cite{Ku91}, raises the question of whether
there exists such a functorial antipode, \ie\ a suitable functor
$\mathbf{S}\colon\mathcal{C}_1\to\mathcal{C}_1$ over $\Vect_\K$.

The existence of such a functor follows from the
equivalence~\eqref{eq_statetkduality} and from the existence of the
homomorphism of Hopf algebras $\underline{S}$ of~\eqref{eq_antipode2},
\ie\ from Proposition~\ref{prop_vertical} applied to
$\mathcal{C}=\comAlg_\K^\op$.
\end{remark}

\begin{proposition}
\label{prop_antipode}
Let
$\mathcal{C}=(\mathcal{C}_0,\mathcal{C}_1,\sigma,\tau,\Uepsilon,\UDelta)$
be a strict symmetric Hopf category according to
Definition~\ref{def_hopfcat}. Here we have used some abbreviations
and written just $\mathcal{C}_0$ for a category
$(\mathcal{C}_0,\omega_0)$ over $\Vect_\K$ and so on. There exists a
functor $\mathbf{S}\colon\mathcal{C}_1\to\mathcal{C}_1$ over
$\Vect_\K$ that satisfies,
\begin{eqnarray}
  \tau\mathbf{S}&=&\sigma,\\
  \sigma\mathbf{S}&=&\tau,\\
  \UDelta{[\mathbf{S};\id_{\mathcal{C}_1}]}_{\id_{\mathcal{C}_0}}\hat\otimes&=&\Uepsilon\sigma,\\
  \UDelta{[\id_{\mathcal{C}_1};\mathbf{S}]}_{\id_{\mathcal{C}_0}}\hat\otimes&=&\Uepsilon\tau,
\end{eqnarray}
where composition of functors is read from left to right.
\end{proposition}

\begin{proof}
Consider the commutative cotrialgebra
$H=\Cat(\coend(-))[\mathcal{C}]=(H_0,H_1,\tilde\sigma,\tilde\tau,\imath,\circ)$
for which $H_1=\coend(\mathcal{C}_1)$. By Remark~\ref{rem_vertical2},
there is a homomorphism of Hopf algebras $\US\colon H_1\to H_1$ as
in~\eqref{eq_antipode2}. Then the functor
$\comod(\US)\colon\comod(H_1)\to\comod(H_1)$ over $\Vect_\K$ yields
the desired functor $\mathbf{S}$ upon choosing an isomorphism
$\mathcal{C}_1\cong\comod(H_1)=\comod(\coend(\mathcal{C}_1))$ in
$\mathfrak{C}^{\otimes,\ast,s}_{\K,rec}$. By Lemma~2.1.3
of~\cite{Sc92}, such an isomorphism is just an equivalence of
categories.
\end{proof}

In order to understand the functorial antipode better than just from
this abstract existence argument, we sketch the situation for the
Hopf category of representations of a $2$-group in
Section~\ref{sect_individual} below.

\subsection{Semisimplicity}

Consider a strict symmetric Hopf category
$\mathcal{C}=(\mathcal{C}_0,\mathcal{C}_1,\sigma,\tau,\Uepsilon,\UDelta)$.
The category $\mathcal{C}_1$ over $\Vect_\K$ is what other authors
would call the Hopf category (Remark~\ref{rem_hopfcat}). There is the
following notion of semisimplicity for symmetric Hopf categories.

\begin{definition}
\begin{myenumerate}
\item
  A $\K$-linear Abelian category $\mathcal{C}$ is called
  \emph{semisimple} if there is a set $\mathcal{C}_0$ of objects of
  $\mathcal{C}$ such that,
\begin{myenumerate}
\item
  any object $X\in\mathcal{C}_0$ is \emph{simple}, \ie\ any non-zero
  monomorphism $f\colon Y\to X$ is an isomorphism and any non-zero
  epimorphism $g\colon X\to W$ is an isomorphism,
\item
  any object $Z$ of $\mathcal{C}$ is isomorphic to an object of the
  form,
\begin{equation}
  Z = \bigoplus_{j=1}^n X_j,
\end{equation}  
  where $n\in\N$ and $X_j\in\mathcal{C}_0$ for all $j$,
\item
  any object $X\in\mathcal{C}_0$ satisfies $\dim_\K\Hom(X,X)=1$.
\end{myenumerate}
  The category $\mathcal{C}$ is called \emph{finitely semisimple} (or
  \emph{Artinian semisimple}) if the set $\mathcal{C}_0$ is finite.
\item
  A reconstructible category $(\mathcal{C},\omega)$ over $\Vect_\K$ is
  called \emph{[finitely] semisimple} if the $\K$-linear Abelian
  category $\mathcal{C}$ is [finitely] semisimple.
\item
  A strict symmetric Hopf category
  $\mathcal{C}=(\mathcal{C}_0,\mathcal{C}_1,\sigma,\tau,\Uepsilon,\UDelta)$
  is called \emph{[finitely]} \emph{semisimple} if both
  $\mathcal{C}_0$ and $\mathcal{C}_1$ are [finitely] semisimple.
\end{myenumerate}
\end{definition}

\begin{definition}
A \emph{$2$-vector space} of Kapranov--Voevodsky type\footnote{A
$2$-vector space of Baez--Crans type~\cite{BaCr04} is an internal
category in $\Vect_\K$.} is a semisimple $\K$-linear Abelian
category. A \emph{finite-dimensional $2$-vector space} is a finitely
semisimple $\K$-linear Abelian category.
\end{definition}

\begin{corollary}
Let
$\mathcal{C}=(\mathcal{C}_0,\mathcal{C}_1,\sigma,\tau,\Uepsilon,\UDelta)$
be a [finitely] semisimple strict symmetric Hopf category. Then the
category $\mathcal{C}$ in $\mathcal{C}_1=(\mathcal{C},\omega)$ is a
[finite-dimensional] $2$-vector space of Kapranov--Voevodsky type.
\end{corollary}

Semisimple Hopf categories originate from cosemisimple cotrialgebras
as follows. For the corepresentation theory underlying the following
results, we refer to~\cite{KlSc97}.

\begin{definition}
\begin{myenumerate}
\item
  A coalgebra $C$ is called \emph{cosemisimple} if it is a direct sum
  of cosimple coalgebras, \ie\ coalgebras that have no non-trivial
  subcoalgebras.
\item
  A strict commutative cotrialgebra
  $H=(H_0,H_1,\sigma,\tau,\Uepsilon,\UDelta)$ is called
  \emph{cosemisimple} if both Hopf algebras $H_0$ and $H_1$ have
  cosemisimple underlying coalgebras.
\end{myenumerate}
\end{definition}

\begin{proposition}[see~\cite{KlSc97}]
Let $C$ be a [finite-dimensional] cosemisimple coalgebra over
$\K$. Then the category $\mathcal{M}^H$ of finite-dimensional
right-$C$-comodules is [finitely] semisimple.
\end{proposition}

\begin{corollary}
Let $H$ be a strict commutative cotrialgebra and
$\mathcal{C}=\Cat(\comod(-))[H]$ be the strict symmetric Hopf category
that is Tannaka--Kre\v\i n dual to $H$. If $H$ is [finite-dimensional
and] cosemisimple as a commutative cotrialgebra, then $\mathcal{C}$ is
[finitely] semisimple as a symmetric Hopf category.
\end{corollary}

\subsection{Representations of compact topological $2$-groups}

The representation theory of a compact topological group $G$ can be
developed as follows by first passing to its commutative Hopf
$C^\ast$-algebra of continuous complex valued functions $C(G)$ and
then studying the comodules of the dense subcoalgebra of
representative functions.

\begin{proposition}[see~\cite{KlSc97}]
Let $G$ be a compact topological group. 
\begin{myenumerate}
\item
  Let $C_\mathrm{alg}(G)$ denote the algebra of representative functions of
  $G$, \ie\ the $\C$-linear span of all matrix elements of
  finite-dimensional continuous unitary representations of $G$. Then
  $C_\mathrm{alg}(G)\subseteq C(G)$ forms a dense subcoalgebra which is
  cosemisimple.
\item
  Each finite-dimensional continuous unitary representation of $G$
  gives rise to a finite-dimensional
  right-$C_\mathrm{alg}(G)$-comodule structure on the same underlying
  complex vector space. Any intertwiner of such representations forms
  a morphism of right-$C_\mathrm{alg}(G)$-comodules.
\item
  Denote by $\hat G$ a set containing one representative for each
  equivalence class of isomorphic simple objects of
  $\mathcal{M}^{C_\mathrm{alg}(G)}$. Then $C_\mathrm{alg}(G)$ is
  isomorphic as a coalgebra to the following direct sum of cosimple
  coalgebras,
\begin{equation}
  C_\mathrm{alg}(G) = \bigoplus_{\rho\in\hat G} V_\rho^\ast\otimes V_\rho,
\end{equation}
  where $V_\rho=\omega^{C_\mathrm{alg}(G)}(\rho)$ denotes the
  underlying vector space of the simple object $\rho$. The vector
  spaces $V_\rho^\ast\otimes V_\rho$ carry the coalgebra structure of
  the cosimple coefficient coalgebras of the simple objects. In
  particular, $C_\mathrm{alg}(G)$ is cosemisimple.
\end{myenumerate}
\end{proposition}

These results show how to proceed in order to formulate the
representation theory of a strict compact topological $2$-group
$G=(G_0,G_1,s,t,\imath,\circ)$. First, use the equivalence of
categories~\eqref{eq_equivcst} in order to obtain the strict
commutative $C^\ast$-cotrialgebra,
\begin{equation}
  H:=\Cat(C(-))[G]=(C(G_0),C(G_1),C(s),C(t),C(\imath),C(\circ)).
\end{equation}
Then restrict to the dense subcoalgebras $C_\mathrm{alg}(G_0)\subseteq
C(G_0)$, \etc. This also restricts the completed tensor products to
the algebraic ones. Notice that this restriction preserves colimits,
\ie\
$H_\mathrm{alg}:=(C_\mathrm{alg}(G_0),C_\mathrm{alg}(G_1),C_\mathrm{alg}(s),C_\mathrm{alg}(t),C_\mathrm{alg}(\imath),C_\mathrm{alg}(\circ))$
is a strict commutative cotrialgebra. Then use the
equivalence~\eqref{eq_twotkduality} in order to obtain a strict
symmetric Hopf category,
\begin{eqnarray}
\label{eq_hopfcatrep}
  \mathcal{C}&:=&\Cat(\comod(-))[H_\mathrm{alg}]\nn\\
  &=& (\comod(C_\mathrm{alg}(G_0)),\comod(C_\mathrm{alg}(G_1)),
     \comod(C_\mathrm{alg}(s)),\comod(C_\mathrm{alg}(t)),\nn\\
  &&\qquad   \comod(C_\mathrm{alg}(\imath)),\comod(C_\mathrm{alg}(\circ))).
\end{eqnarray}
This strict symmetric Hopf category $\mathcal{C}$ plays the role of
the category of `finite-dimensional continuous unitary
representations' of the compact topological $2$-group. Recall that
$C_\mathrm{alg}(-)$ is a covariant functor if written as
$C_\mathrm{alg}(-)\colon\compTopGrp\to\comHopf_\C^\op$, and so
$\mathcal{C}$ is indeed an internal cocategory in
$\mathfrak{C}^{\otimes,\ast,s}_{\C,rec}$.

\subsection{Individual representations of $2$-groups}
\label{sect_individual}

Let $G$ be a compact topological group and $\Rep
(G):=\comod(C_\mathrm{alg}(G))$ be its category of finite-dimensional
continuous unitary representations, viewed as
right-$C_\mathrm{alg}(G)$-comodules with the forgetful functor to
$\fdVect_\C$.

If $f\colon G\to H$ is a homomorphism of compact topological groups,
then $\Rep(f)\colon\Rep(H)\to\Rep(G)$ is the functor that assigns to each
finite-dimensional continuous representation of $H$ the same
underlying vector space, but viewed via $f$ as a representation of
$G$.

The strict symmetric Hopf category~\eqref{eq_hopfcatrep} is in this
notation,
\begin{eqnarray}
  \Rep(G) &=& (\Rep(G_0),\Rep(G_1),\Rep(s),\Rep(t),\Rep(\imath),\Rep(\circ))\nn\\
  &:=&\Cat(\comod (-))[H_\mathrm{alg}].
\end{eqnarray}
It plays the role of the representation \emph{category} of the strict
compact topological $2$-group $G=(G_0,G_1,s,t,\imath,\circ)$. What is
an \emph{individual} representation of $G$?

The answer is that an individual representation is just a
(finite-dimensional continuous unitary) representation of $G_1$. There
is no difference compared with the representations of ordinary groups
except that for $2$-groups, the representation category $\Rep(G_1)$
has got more structure. The purpose of the category $\Rep(G_0)$ and of
the functors $\Rep(s)$, $\Rep(t)$, $\Uepsilon:=\Rep(\imath)$ and
$\UDelta:=\Rep(\circ)$ is merely to describe this additional
structure.

The most important additional structure is the functorial
comultiplication, \ie\ the symmetric monoidal functor
\begin{equation}
  \UDelta\colon\Rep(G_1)\to\Rep(G_1)\pushout{\Rep(s)}{\Rep(G_0)}{\Rep(t)}\Rep(G_1)
    \cong\Rep\bigl(G_1\pullback{s}{G_0}{t}G_1\bigr),
\end{equation}
where the isomorphism sign `$\cong$' denotes isomorphism of objects of
$\mathfrak{C}^{\otimes,\ast,s}_{\K,rec}$, \ie\ equivalence of the
corresponding categories by Lemma~2.1.3 of~\cite{Sc92}. The functor
$\UDelta$ over $\Vect_\K$ assigns to each finite-dimensional unitary
continuous representation of $G_1$ the same underlying vector space,
but viewed as representation of $G_1\pullback{s}{G_0}{t}G_1$ via the
homomorphism
\begin{equation}
  \circ\colon G_1\pullback{s}{G_0}{t}G_1=\{\,(g,g^\prime)\in G_1\times G_1\colon\quad
    s(g)=t(g^\prime)\,\}\to G_1.
\end{equation}
A functorial antipode can be constructed as in
Proposition~\ref{prop_antipode}. It is the functor
\begin{equation}
  \Rep(\xi)\colon\Rep(G_1)\to\Rep(G_1)
\end{equation}
over $\Vect_\K$ which is obtained from vertical inversion
(Proposition~\ref{prop_vertical}(2)) and which is specified up to natural
isomorphism over $\Vect_\K$.

%
\section{The non-commutative and non-symmetric cases}
%
\label{sect_gentrialg}

Given the construction of cocommutative trialgebras and commutative
cotrialgebras in Sections~\ref{sect_coctrialg}
and~\ref{sect_comcotrialg}, one might be tempted to try the following
definition.

\begin{definition}
\label{def_generic}
The objects of the $2$-category,
\begin{equation}
  \Tri_\K := \Cat(\Hopf_\K),
\end{equation}
are called \emph{strict generic trialgebras} and the objects of,
\begin{equation}
  \CoTri_\K := \CoCat(\Hopf_\K),
\end{equation}
\emph{strict generic cotrialgebras}.
\end{definition}

\begin{remark}
The above definition comes with an important warning, though. Although
each cocommutative Hopf algebra is at the same time a Hopf algebra, a
cocommutative trialgebra would not necessarily be a generic
trialgebra. This is a consequence of the fact that the inclusion
functor $\cocHopf_\K\to\Hopf_\K$ does not preserve all finite
limits. Similarly, a commutative cotrialgebra would not necessarily be
a generic cotrialgebra. In particular, strict $2$-groups and crossed
modules are not immediately helpful in constructing examples of
generic (co-)trialgebras.

A second objection against Definition~\ref{def_generic} is the
observation that it would use the categorical product and pullback in
$\Hopf_\K$ which is in general `much bigger' than the tensor product
of Hopf algebras. Experience with the theory of Hopf algebras,
however, suggests that one ought to consider their tensor product
rather than their categorical product.

Notice further that in general, although every Hopf algebra has got an
underlying vector space, neither the objects of $\Cat(\cocHopf_\K)$
nor the objects of $\CoCat(\comHopf_\K)$ come with an underlying
strict $2$-vector space of Baez--Crans type~\cite{BaCr04}. A strict
$2$-vector space of that sort is an object of $\Cat(\Vect_\K)$. This
is a consequence of the fact that the forgetful functors
$\cocHopf_\K\to\Vect_\K$ and $\comHopf_\K\to\Vect_\K$ do not preserve
all finite limits.
\end{remark}

As soon as one has an example of a strict generic cotrialgebra,
Tannaka--Kre\v\i n duality can still be used in order to obtain strict
generic Hopf categories. In order to show this, we drop the
requirement of symmetry from Definition~\ref{def_catvect} and proceed
as follows.

\begin{definition}
\begin{myenumerate}
\item
  $\mathfrak{C}^{\otimes,\ast}_{\K,rec}$ denotes the category whose
  objects are rigid monoidal categories $(\mathcal{C},\omega)$ over
  $\Vect_\K$ such that $\mathcal{C}$ is a $\K$-linear Abelian
  essentially small category and $\omega\colon\mathcal{C}\to\Vect_\K$
  is an exact faithful $\K$-linear functor into $\fdVect_\K$. The
  morphisms of $\mathfrak{C}^{\otimes,\ast}_{\K,rec}$ are monoidal
  functors over $\Vect_\K$.
\item
  The objects of the $2$-category
\begin{equation}
  \HopfCat_\K := \CoCat(\mathfrak{C}^{\otimes,\ast}_{\K,rec})
\end{equation}
  are called \emph{strict generic Hopf categories over} $\K$.
\end{myenumerate}
\end{definition}
Then, slightly modifying the proofs of Appendix~\ref{app_tkduality},
we obtain Tannaka--Kre\v\i n duality in the following form.

\begin{theorem}
There is an equivalence of categories,
\begin{equation}
  \mathfrak{C}^{\otimes,\ast}_{\K,rec}\simeq\Hopf_\K,
\end{equation}
given by the functors
\begin{equation}
  \comod(-)\colon\Hopf_\K\to\mathfrak{C}^{\otimes,\ast}_{\K,rec}
\end{equation}
and
\begin{equation}
  \coend(-)\colon\mathfrak{C}^{\otimes,\ast}_{\K,rec}\to\Hopf_\K.
\end{equation}
\end{theorem}

\begin{corollary}
There is a $2$-equivalence between the $2$-categories
\begin{equation}
  \CoTri_\K\simeq\HopfCat_\K
\end{equation}
provided by the functors
\begin{equation}
  \Cat(\comod(-))\colon\CoTri_\K\to\HopfCat_\K
\end{equation}
and
\begin{equation}
  \Cat(\coend(-))\colon\HopfCat_\K\to\CoTri_\K.
\end{equation}
\end{corollary}

In order to generalize both the notions of trialgebra and of Hopf
categories beyond the strict case, one strategy would be to set up the
entire analysis of the present article for the more general case of
weak or coherent $2$-groups~\cite{BaLa04}. So far, it is not obvious
whether these weaker structures form models of any algebraic or
essentially algebraic theory in $\Grp$, and so there is no obvious way
of employing functors such as $k[-]\colon\Grp\to\cocHopf_\K$. It seems
that one either needs to extend and generalize the relevant universal
algebra first or that one has to find a different way of relating the
various algebraic structures.

\acknowledgments

I would like to thank John Baez, Aaron Lauda and Karl-Georg
Schlesinger for valuable discussions and the Erwin Schr\"odinger
Institute, Vienna, for hosting the 2004 Programme on Tensor Categories
where this project was started.

%
\appendix
%

\section{Compact topological groups}
\label{app_topgrp}

In this appendix, we sketch how to encapsulate all the functional
analysis of Gel'fand representation theory in the algebraic statement
that the category $\compHaus$ of compact Hausdorff spaces with
continuous maps is equivalent to the opposite of the category
$\comUnCstAlg$ of commutative unital $C^\ast$-algebras with unital
$\ast$-homomorphisms,
\begin{equation}
  \compHaus\simeq\comUnCstAlg^\op.
\end{equation}
This allows us to present concise definitions for the category
$\compTopGrp$ of compact topological groups and for the category
$\comHopfCstAlg$ of commutative Hopf $C^\ast$-algebras so that we can
establish the equivalence,
\begin{equation}
  \compTopGrp\simeq\comHopfCstAlg^\op,
\end{equation}
which we have stated as Theorem~\ref{thm_usegelfand}.

\subsection{Gel'fand representation theory}

For background on $C^\ast$-algebras and for the proofs of the results
summarized here, see, for example~\cite{Mu90}.

There is a category $\comUnCstAlg$ whose objects are commutative
unital $C^\ast$-algebras and whose morphisms are unital
$\ast$-homomorphisms. This category has all finite coproducts. In
particular, its initial object is the field of complex numbers
$\C$. For each commutative unital $C^\ast$-algebra $A$, there is a
unique unital $\ast$-homomorphism $\C\to A$, $1_\C\mapsto 1_A$. The
binary coproduct $A\coprod B=A\otimes_\ast B$ of two commutative
unital $C^\ast$-algebras $A$ and $B$ is the completion of the tensor
product in the so-called spatial $C^\ast$-norm ${||\cdot||}_\ast$. By
a theorem of Takesaki, this is the unique $C^\ast$-norm on the
algebraic tensor product $A\otimes B$ whereby the uniqueness result
exploits that $A$ and $B$ are commutative algebras. The colimiting
cone of the coproduct is given by the unital $\ast$-homomorphisms
$\imath_A\colon A\to A\otimes_\ast B$, $a\mapsto a\otimes 1_B$ and
$\imath_B\colon B\to A\otimes_\ast B$, $b\mapsto 1_A\otimes b$.

\begin{proposition}
There is a functor $C(-)\colon\compHaus\to\comUnCstAlg^\op$ which
assigns to each compact Hausdorff space $X$ the commutative unital
$C^\ast$-algebra $C(X)$ of continuous complex valued functions on $X$
with the supremum norm,
\begin{equation}
  ||f|| := \sup_{x\in X} |f(x)|,\quad f\in C(X).
\end{equation}
The functor $C(-)$ assigns to each continuous map $\phi\colon X\to Y$
between compact Hausdorff spaces the unital $\ast$-homomorphism,
\begin{equation}
  C(\phi)\colon C(Y)\to C(X),f\mapsto \phi f
\end{equation}
(first $\phi$, then $f$).
\end{proposition}

Let $A$ be a commutative $C^\ast$-algebra. A \emph{character} of $A$
is a non-zero linear functional $\omega\colon A\to\C$ for which
$\omega(ab)=\omega(a)\omega(b)$ for all $a,b\in A$. One can show that
characters are continuous maps and that they satisfy $\omega(1_A)=1$
if $A$ is unital. The \emph{Gel'fand spectrum} $\sigma(A)$ of $A$ is
the set of all characters of $A$. For each $a\in A$, one defines the
Gel'fand transform $\hat a\colon\sigma(A)\to\C$, $\omega\to
\omega(a)$.

\begin{proposition}
There is a functor $\sigma(-)\colon\comUnCstAlg\to\compHaus^\op$ which
assigns to each commutative unital $C^\ast$-algebra $A$ its Gel'fand
spectrum $\sigma(A)$. The set $\sigma(A)$ forms a compact Hausdorff
space if it is equipped with the weak${}^\ast$ topology which is the
coarsest topology for $\sigma(A)$ such that all maps $\hat
a\colon\sigma(A)\to\C$, $a\in A$, are continuous. The functor
$\sigma(-)$ assigns to each unital $\ast$-homomorphism $\Phi\colon
A\to B$ between commutative unital $C^\ast$-algebras $A$ and $B$ the
continuous map,
\begin{equation}
  \sigma(\Phi)\colon\sigma(B)\to\sigma(A), \omega\mapsto\Phi\omega
\end{equation}
(first $\Phi$, then $\omega$).
\end{proposition}

\begin{theorem}[Gel'fand]
\label{thm_gelfand}
Let $A$ be a commutative unital $C^\ast$-algebra. The Gel'fand
transform,
\begin{equation}
  \hat{\ }\colon A\to C(\sigma(A)), a\mapsto \hat a,
\end{equation}
is a unital $\ast$-isomorphism.
\end{theorem}

This theorem shows that each commutative unital $C^\ast$-algebra
arises as the algebra of functions on its Gel'fand spectrum. In order
to formulate the complementary result that each compact Hausdorff
space $X$ arises as the Gel'fand spectrum of its function algebra, we
define the following characters of $C(X)$ by evaluation at $x\in X$,
\begin{equation}
  \omega_x\colon C(X)\to\C, f\mapsto f(x).
\end{equation}

\begin{theorem}[Gel'fand]
\label{thm_gelfandinv}
Let $X$ be a compact Hausdorff space. The map,
\begin{equation}
  \omega\colon X\to\sigma(C(X)), x\mapsto \omega_x,
\end{equation}
is a homeomorphism of topological spaces.
\end{theorem}

The analytical results reviewed in this section can be summarized in
categorical terms as follows.

\begin{theorem}
\label{thm_gelfandall}
There is an equivalence of categories
$\compHaus\simeq\comUnCstAlg^\op$ provided by the functors
\begin{equation}
  C(-)\colon\compHaus\to\comUnCstAlg^\op
\end{equation}
and
\begin{equation}
  \sigma(-)\colon\comUnCstAlg^\op\to\compHaus
\end{equation}
with the natural isomorphisms $\hat{\ }\colon
1_{\comUnCstAlg^\op}\Rightarrow \sigma(-)\C(-)$ defined in
Theorem~\ref{thm_gelfand} and $\omega\colon 1_\compHaus\Rightarrow
C(-)\sigma(-)$ defined in Theorem~\ref{thm_gelfandinv}.
\end{theorem}

Recall that the category $\compHaus$ has all finite products. In
particular, the terminal object is the one-element topological space
$\{\ast\}$ and the binary product $X\prod Y=X\times Y$ of two compact
Hausdorff spaces $X$ and $Y$ is their Cartesian product with the
projection maps $p_1\colon X\times Y\to X, (x,y)\mapsto x$ and
$p_2\colon X\times Y\to Y, (x,y)\mapsto y$. We can now use
Theorem~\ref{thm_gelfandall} in order to relate products in
$\compHaus$ with coproducts in $\comUnCstAlg$. In particular there are
isomorphisms of unital $\ast$-algebras,
\begin{equation}
  C(\{\ast\})\cong \C\qquad\mbox{and}\qquad 
  C(X\times Y)\cong C(X)\otimes_\ast C(Y)
\end{equation}
and homeomorphisms of topological spaces,
\begin{equation}
  \sigma(\C)\cong \{\ast\}\qquad\mbox{and}\qquad
  \sigma(A\otimes_\ast B)\cong \sigma(A)\times\sigma(B).
\end{equation}

\subsection{Group objects}
\label{app_gpobj}

In this section, we recall the definition of group objects in
categories with finite products. Then we can define the notion of a
topological group as a group object in $\compHaus$ and of a
commutative Hopf $C^\ast$-algebra as a group object in
$\comUnCstAlg^\op$.

\begin{definition}
\label{def_groupobject}
Let $\mathcal{C}$ be a category with finite products and $T$ be a
terminal object of $\mathcal{C}$.
\begin{myenumerate}
\item
  A \emph{group object} $G=(C,\mu,\eta,\zeta)$ in $\mathcal{C}$
  consists of an object $C$ of $\mathcal{C}$ and of morphisms
  $\mu\colon C\prod C\to C$, $\eta\colon T\to C$ and $\zeta\colon C\to
  C$ of $\mathcal{C}$ such that the following diagrams commute,
\begin{equation}
\label{eq_groupobj1}
\begin{aligned}
\xymatrix{
  C\prod C\prod C\ar[rr]^{(\mu,\id_C)}\ar[dd]_{(\id_C,\mu)}&&C\prod C\ar[dd]^\mu\\
  \\
  C\prod C\ar[rr]_\mu&&C
}
\end{aligned}
\end{equation}
\begin{equation}
\label{eq_groupobj2}
\begin{aligned}
\xymatrix{
  T\prod C\ar[rr]^{(\eta,\id_C)}\ar[rrdd]_{p_2}&&
    C\prod C\ar[dd]^\mu&&C\prod T\ar[ll]_{(\id_C,\eta)}\ar[lldd]^{p_1}\\
  \\
  &&C
}
\end{aligned}
\end{equation}
\begin{equation}
\label{eq_groupobj3}
\begin{aligned}
\xymatrix{
  C\ar[rr]\ar[dd]_{\delta(\id_C,\zeta)}&&T\ar[dd]^\eta\\
  \\
  C\prod C\ar[rr]_\mu&& C
}
\end{aligned}\qquad\qquad
\begin{aligned}
\xymatrix{
  C\ar[rr]\ar[dd]_{\delta(\zeta,\id_C)}&&T\ar[dd]^\eta\\
  \\
  C\prod C\ar[rr]_\mu&& C
}
\end{aligned}
\end{equation}
\item
  Let $G=(C,\mu,\eta,\zeta)$ and
  $G^\prime=(C^\prime,\mu^\prime,\eta^\prime,\zeta^\prime)$ be group
  objects in $\mathcal{C}$. An \emph{internal group homomorphism}
  $\phi\colon G\to G^\prime$ is a morphism $\phi\colon C\to C^\prime$
  of $\mathcal{C}$ such that the following diagrams commute,
\begin{equation}
\label{eq_intgrphom1}
\begin{aligned}
\xymatrix{
  C\prod C\ar[rr]^\mu\ar[dd]_{(\phi,\phi)}&&C\ar[dd]^\phi\\
  \\
  C^\prime\prod C^\prime\ar[rr]_{\mu^\prime}&&C^\prime
}
\end{aligned}
\end{equation}
\begin{equation}
\label{eq_intgrphom2}
\begin{aligned}
\xymatrix{
  T\ar[rr]^\eta\ar[rrdd]_{\eta^\prime}&&C\ar[dd]^\phi\\
  \\
  &&C^\prime  
}
\end{aligned}
\end{equation}
\begin{equation}
\label{eq_intgrphom3}
\begin{aligned}
\xymatrix{
  C\ar[rr]^\zeta\ar[dd]_\phi&&C\ar[dd]^\phi\\
  \\
  C^\prime\ar[rr]_{\zeta^\prime}&&C^\prime
}
\end{aligned}
\end{equation}
\item
  Let $G=(C,\mu,\eta,\zeta)$,
  $G^\prime=(C^\prime,\mu^\prime,\eta^\prime,\zeta^\prime)$ and
  $G^{\prime\prime}=(C^{\prime\prime},\mu^{\prime\prime},\eta^{\prime\prime},
  \zeta^{\prime\prime})$ be group objects in $\mathcal{C}$ and
  $\phi\colon G\to G^\prime$ and $\psi\colon G^\prime\to
  G^{\prime\prime}$ be internal group homomorphisms in
  $\mathcal{C}$. The \emph{composition} of $\phi$ and $\psi$ is the
  internal group homomorphism $\phi\cdot\psi\colon G\to
  G^{\prime\prime}$ defined by the morphism $\phi\cdot\psi\colon C\to
  C^{\prime\prime}$ of $\mathcal{C}$.
\item
  Let $G=(C,\mu,\eta,\zeta)$ be a group object in $\mathcal{C}$. The
  \emph{identity internal group homomorphism} $\id_G\colon G\to G$ is
  defined as the morphism $\id_C\colon C\to C$ of $\mathcal{C}$.
\end{myenumerate}
\end{definition}

\begin{theorem}
Let $\mathcal{C}$ be a category with finite products. There is a
category $\Grp(\mathcal{C})$ whose objects are group objects and whose
morphisms are internal group homomorphisms in $\mathcal{C}$.
\end{theorem}

Similarly to the definition of internal categories in
Section~\ref{sect_internalization}, the guiding example is the
following.

\begin{example}
The category $\Set$ has finite products. $\Grp(\Set)$ is the category
of groups and group homomorphisms.
\end{example}

\begin{remark}
\label{rem_technical}
In the definition of group objects and internal group homomorphisms
(Definition~\ref{def_groupobject}) and also in the definition of
internal categories, functors and natural transformations
(Definition~\ref{def_internalcat}), there is a technical subtlety
which we have suppressed.

All objects that are constructed as limits and which are used in the
diagrams~\eqref{eq_groupobj1} to~\eqref{eq_intgrphom3}, for example,
the terminal object $T$ or the product $C\prod C$, are defined only up
to isomorphism. The definition of a group object in $\mathcal{C}$
should therefore contain an object $C$ of $\mathcal{C}$ \emph{and in
addition} the choice of an object $T$ of $\mathcal{C}$ which is
terminal and of an object $C^2$ isomorphic to the product $C\prod C$,
and so on. Correspondingly, the definition of an internal group
homomorphism should contain besides the morphism $\phi\colon C\to
C^\prime$ \emph{in addition} a morphism $\phi^2\colon C^2\to C^2$, and
so on. This is tidied up and properly taken into account, for example,
by the construction of \emph{sketches} as in~\cite{BaWe83}. We have
been reluctant to do the same for internal categories in
Section~\ref{sect_internalization} since this would have made that
section far less accessible.
\end{remark}

Similarly to the study of internal categories in
Section~\ref{sect_internalization}, we are interested in varying the
base category $\mathcal{C}$ in the definition of $\Grp(\mathcal{C})$.

\begin{proposition}
Let $\mathcal{C}$ and $\mathcal{D}$ be categories with finite
products and $T\colon\mathcal{C}\to\mathcal{D}$ be a finite-product
preserving functor. Then there is a functor
$\Grp(T)\colon\Grp(\mathcal{C})\to\Grp(\mathcal{D})$ given as follows.
\begin{myenumerate}
\item
  $\Grp(T)$ associates with each group object $G=(C,\mu,\eta,\zeta)$
  in $\mathcal{C}$ the group object
  $\Grp(T)[G]:=(TC,T\mu,T\eta,T\zeta)$ in $\mathcal{D}$.
\item
  Let $G=(C,\mu,\eta,\zeta)$ and
  $G^\prime=(C^\prime,\mu^\prime,\eta^\prime,\zeta^\prime)$ be group
  objects in $\mathcal{C}$. $\Grp(T)$ associates with each internal
  group homomorphism $\phi\colon G\to G^\prime$ in $\mathcal{C}$ the
  internal group homomorphism
  $\Grp(T)[\phi]\colon\Grp(\mathcal{C})\to\Grp(\mathcal{D})$ in
  $\mathcal{D}$ which is given by the morphism $T\phi\colon TC\to
  TC^\prime$ of $\mathcal{D}$.
\end{myenumerate}
\end{proposition}

\begin{proposition}
Let $\mathcal{C}$ and $\mathcal{D}$ be categories with finite
products, $T,\tilde T\colon\mathcal{C}\to\mathcal{D}$ be
finite-product preserving functors and $\alpha\colon
T\Rightarrow\tilde T$ a natural transformation. Then there is a
natural transformation
$\Grp(\alpha)\colon\Grp(T)\Rightarrow\Grp(\tilde T)$. It associates
with each group object $G=(C,\mu,\eta,\zeta)$ in $\mathcal{C}$ the
internal group homomorphism
${\Grp(\alpha)}_G\colon\Grp(T)[G]\to\Grp(\tilde T)[G]$ in
$\mathcal{D}$ which is given by the morphism $\alpha_C\colon
TC\to\tilde TC$.
\end{proposition}

\begin{theorem}
Let $\fpCat$ denote the $2$-category of small categories with finite
products, finite-product preserving functors and natural
transformations. Let $\Cat$ denote the $2$-category of small
categories, functors and natural transformations. Then $\Grp(-)$ forms
a $2$-functor,
\begin{equation}
  \Grp(-)\colon\fpCat\to\Cat.
\end{equation}
\end{theorem}

\begin{remark}
It is known that the functor $\Grp(-)$ is actually a functor
$\Grp(-)\colon\fpCat\to\fpCat$ and even
$\Grp(-)\colon\fcCat\to\fcCat$. In particular, if $\mathcal{C}$ has
all finite products, so does the category $\Grp(\mathcal{C})$, and if
$\mathcal{C}$ has all finite limits, so does
$\Grp(\mathcal{C})$~\cite{BaWe83}.
\end{remark}

\begin{corollary}
\label{corr_groupobj}
Let $\mathcal{C}\simeq\mathcal{D}$ be an equivalence of categories
with finite products provided by (finite-product preserving) functors
$F\colon\mathcal{C}\to\mathcal{D}$ and
$G\colon\mathcal{D}\to\mathcal{C}$ with natural isomorphisms
$\eta\colon 1_\mathcal{C}\Rightarrow FG$ and $\epsilon\colon
GF\Rightarrow 1_\mathcal{D}$. Then there is an equivalence of
categories $\Grp(\mathcal{C})\simeq\Grp(\mathcal{D})$ given by the
functors $\Grp(F)\colon\Grp(\mathcal{C})\to\Grp(\mathcal{D})$ and
$\Grp(G)\colon\Grp(\mathcal{D})\to\Grp(\mathcal{C})$ with natural
isomorphisms $\Grp(\eta)\colon 1_{\Grp(\mathcal{C})}\Rightarrow
\Grp(F)\Grp(G)$ and $\Grp(\epsilon)\colon\Grp(G)\Grp(F)\Rightarrow
1_{\Grp(\mathcal{D})}$.
\end{corollary}

\begin{definition}
Let $\mathcal{C}$ be a category with finite coproducts. We define,
\begin{equation}
\label{eq_cogrp}
  \CoGrp(\mathcal{C}) := {\Grp(\mathcal{C}^\op)}^\op,
\end{equation}
and call the objects of $\CoGrp(\mathcal{C})$ \emph{cogroup objects}
in $\mathcal{C}$ and the morphisms \emph{internal cogroup
homomorphisms} in $\mathcal{C}$.
\end{definition}

\begin{remark}
  We have added the last `$\op$' in~\eqref{eq_cogrp} in order to make
  some familiar functors turn out to be contravariant as they are
  usually written. Note that such an `$\op$' is not present in the
  analogous definition of an internal cocategory in
  Definition~\ref{def_cocat}.
\end{remark}

\begin{definition}
\label{def_hopfcstar2}
The objects of the category,
\begin{equation}
  \compTopGrp := \Grp(\compHaus),
\end{equation}
are called \emph{compact topological groups}. The objects of,
\begin{equation}
  \comHopfCstAlg := \CoGrp(\comUnCstAlg),
\end{equation}
are called \emph{commutative Hopf $C^\ast$-algebras}.
\end{definition}

\begin{remark}
\begin{myenumerate}
\item
  Let $f,g$ be morphisms of $\compHaus$, then $(f;g)=f\times g$ (\cf\
  Section~\ref{sect_notation}). For morphisms $k,\ell$ of
  $\comUnCstAlg$, $[k;\ell]=k\otimes_\ast\ell$. Let $A$ be an object
  of $\comUnCstAlg$, then $\delta^\op\colon A\otimes_\ast A\to A$ is
  precisely the multiplication operation of $A$ since the unit law of
  the multiplication agrees with the defining condition of
  $\delta^\op$.
\item
  Comparing Definition~\ref{def_hopfcstar} with
  Definition~\ref{def_hopfcstar2}, we note that~\eqref{eq_groupobj1}
  gives the coassociativity axiom, \eqref{eq_groupobj2} gives the
  counit axiom and \eqref{eq_groupobj3} gives the antipode axiom.
\end{myenumerate}
\end{remark}

\subsection{Gel'fand representation theory for compact topological groups}

Combining Gel'fand representation theory
(Theorem~\ref{thm_gelfandall}) with Corollary~\ref{corr_groupobj}, we
obtain the main result of this Appendix:

\begin{theorem}
There is an equivalence of categories,
\begin{equation}
  \compTopGrp=\Grp(\compHaus)\simeq\Grp(\comUnCstAlg^\op)=\comHopfCstAlg^\op
\end{equation}
provided by the functors,
\begin{equation}
  \Grp(C(-))\colon\compTopGrp\to\comHopfCstAlg^\op,
\end{equation}
and
\begin{equation}
  \Grp(\sigma(-))\colon\comHopfCstAlg^\op\to\compTopGrp.
\end{equation}
\end{theorem}

In the formulation as Theorem~\ref{thm_usegelfand}, we have just
omitted the $\Grp(-)$ in the name of the functors in order to keep the
notation simple.

\section{Tannaka--Kre\v\i n duality}
\label{app_tkduality}

In this appendix, we review the key results on Tannaka--Kre\v\i n
duality following the presentation of Schauenburg~\cite{Sc92} and show
how to derive Theorem~\ref{thm_statetkduality}. We refer heavily
to~\cite{Sc92} in order to take the quickest route to the theorem. In
the following, $\K$ is some fixed field.

Schauenburg~\cite{Sc92} shows the existence of an adjunction between
the category of Hopf algebras over $\K$ and a category of monoidal
categories over $\Vect_\K$ which we restrict to the following
category.

\begin{definition}
$\mathfrak{C}^{\otimes,\ast}_{\K,rec}$ denotes the category whose
objects are rigid monoidal categories $(\mathcal{C},\omega)$ over
$\Vect_\K$ such that $\mathcal{C}$ is a $\K$-linear Abelian
essentially small category and $\omega\colon\mathcal{C}\to\Vect_\K$ an
exact faithful $\K$-linear functor with values in $\fdVect_\K$. The
morphisms of $\mathfrak{C}^{\otimes,\ast}_{\K,rec}$ are monoidal
functors over $\Vect_\K$.
\end{definition}

\begin{theorem}[see~\cite{Sc92}]
There is an adjunction,
\begin{equation}
\label{eq_schauadjunction}
\begin{aligned}
\xymatrix{
  \mathfrak{C}^{\otimes,\ast}_{\K,rec}\ar@/^2ex/[rr]^{\coend(-)}="p"&&\Hopf_\K\ar@/^2ex/[ll]^{\comod(-)}="q"
      \ar@{-|}"p"+<0ex,-2.5ex>;"q"+<0ex,2.5ex>
}
\end{aligned}
\end{equation}
\end{theorem}

\begin{proof}
This result is stated in Remark~2.4.4 of~\cite{Sc92} for a category
that is bigger than our $\mathfrak{C}^{\otimes,\ast,s}_{\K,rec}$. For
the definition of the functors, see
Theorem~\ref{thm_statetkduality}. The unit of the adjunction is given
by a functor
$[I_\mathcal{C},\id_\omega]\colon(\mathcal{C},\omega)\to(\mathcal{M}^H,\omega^H)$
over $\Vect_\K$ for each object $(\mathcal{C},\omega)$ of
$\mathfrak{C}^{\otimes,\ast}_{\K,rec}$. Here
$H:=\coend(\mathcal{C},\omega)$ and
$(\mathcal{M}^H,\omega^H):=\comod(H)$. The functor
$I_\mathcal{C}\to\mathcal{M}^H$ is defined in Theorem~2.1.12
of~\cite{Sc92}. The counit is given by bialgebra homomorphisms
$\epsilon_H\colon\coend(\mathcal{M}^H,\omega^H)\to H$ where $H$ is any
Hopf algebra over $\K$ (Lemma~2.2.1 of~\cite{Sc92}). The $\epsilon_H$
are obtained from the universal property of the coendomorphism
coalgebra by diagram completion.

In addition to~\cite{Sc92}, we have restricted the category on the
left hand side of~\eqref{eq_schauadjunction} to our
$\mathfrak{C}^{\otimes,\ast}_{\K,rec}$ by adding the
`reconstructibility' conditions that any object $(\mathcal{C},\omega)$
of $\mathfrak{C}^{\otimes,\ast}_{\K,rec}$ consist of a $\K$-linear
Abelian essentially small category $\mathcal{C}$ and an exact faithful
$\K$-linear functor $\omega$ (Section~2.2 of~\cite{Sc92}). This also
guarantees that the category of monoid objects which appear
in~\cite{Sc92} on the left hand side, indeed agrees with our
definition of $\mathfrak{C}^{\otimes,\ast}_{\K,rec}$, \cf\ Lemma~2.3.4
of~\cite{Sc92} and the comments thereafter.
\end{proof}

We now restrict the adjunction~\eqref{eq_schauadjunction} to the
commutative and symmetric case. Note that
$\mathfrak{C}^{\otimes,\ast,s}_{\K,rec}$ is a full subcategory of
$\mathfrak{C}^{\otimes,\ast}_{\K,rec}$
(Proposition~\ref{prop_symmetric}) and that $\comHopf_\K$ is a full
subcategory of $\Hopf_\K$. In order to retain the adjunction, we have
to show that both functors restrict (on objects) to the subcategories.

\begin{proposition}
\begin{myenumerate}
\item
  Let $H$ be a commutative Hopf algebra over $\K$. Then the category
  $\mathcal{M}^H$ of finite-dimensional right-$H$-comodules forms a
  $\K$-linear Abelian rigid symmetric monoidal category. The forgetful
  functor $\omega^H\colon\mathcal{M}^H\to\fdVect_\K$ is a $\K$-linear
  exact faithful symmetric monoidal functor.
\item
  Let $(\mathcal{C},\omega)$ be an object of
  $\mathfrak{C}^{\otimes,\ast,s}_{\K,rec}$. Then the reconstructed
  Hopf algebra $H:=\coend(\mathcal{C},\omega)$ is commutative.
\end{myenumerate}
\end{proposition}

\begin{proof}
The first part is standard. For the second part, commutativity of $H$
is a consequence of $\omega$ being a \emph{symmetric} monoidal
functor. The proof involves a large commutative diagram.
\end{proof}

\begin{corollary}
\label{corr_symmetric}
There is an adjunction,
\begin{equation}
\label{eq_adjunction}
\begin{aligned}
\xymatrix{
  \mathfrak{C}^{\otimes,\ast,s}_{\K,rec}\ar@/^2ex/[rr]^{\coend(-)}="p"&&\comHopf_\K\ar@/^2ex/[ll]^{\comod(-)}="q"
      \ar@{-|}"p"+<0ex,-2.5ex>;"q"+<0ex,2.5ex>
}
\end{aligned}
\end{equation}
\end{corollary}

We finally show that both unit and counit of the adjunction are
natural isomorphisms and thereby establish the equivalence of
categories claimed in Theorem~\ref{thm_statetkduality}.

\begin{proposition}
The adjunction~\eqref{eq_adjunction} is an equivalence of categories.
\end{proposition}

\begin{proof}
In order to see that the counit is a natural isomorphism, let $H$ be a
commutative Hopf algebra over $\K$. The counit
$\epsilon_H\colon\coend(\mathcal{M}^H,\omega^H)\to H$ forms an
isomorphism of coalgebras (Lemma~2.2.1 of~\cite{Sc92}). Since
$\epsilon_H$ is also a homomorphism of bialgebras, it forms an
isomorphism of bialgebras and therefore also an isomorphism in the
category $\comHopf_\K$.

In order to see that the unit is a natural isomorphism, let
$(\mathcal{C},\omega)$ be an object of
$\mathfrak{C}^{\otimes,\ast,s}_{\K,rec}$. The conditions on
$(\mathcal{C},\omega)$ then guarantee that the functor
$I_\mathcal{C}\to\mathcal{M}^H$, $H=\coend(\mathcal{C},\omega)$ that
features in the unit, forms an equivalence of categories (Section~2.2
of~\cite{Sc92}). Together with the identity natural isomorphism of the
underlying forgetful functor, it forms a functor
$[I_\mathcal{C},\id_\omega]\colon(\mathcal{C},\omega)\to(\mathcal{M}^H,\omega^H)$
over $\Vect_\K$. By Corollary~2.3.7 of~\cite{Sc92}, $I_\mathcal{C}$ is
a monoidal functor. Since $\id_\omega$ is a monoidal natural
transformation, $[I_\mathcal{C},\id_\omega]$ is a monoidal functor
over $\Vect_\K$, and by our Proposition~\ref{prop_symmetric} it is
symmetric monoidal as a functor over $\Vect_\K$.

Consider now the representative $(I_\mathcal{C},\id_\omega)$ of the
equivalence class $[I_\mathcal{C},\id_\omega]$. Since the functor
$\omega^H\colon\mathcal{M}^H\to\Vect_\K$ is the underlying forgetful
functor, the condition~\eqref{eq_symmetricfunctor} that
$[I_\mathcal{C},\id_\omega]$ is a symmetric monoidal functor over
$\Vect_\K$, implies that
$I_\mathcal{C}\colon\mathcal{C}\to\mathcal{M}^H$ is a symmetric
monoidal functor. Since $I_\mathcal{C}$ is part of an equivalence of
categories, by~\cite{Sa72}, Chapter~I, Proposition~4.4.2,
$I_\mathcal{C}$ forms a tensor equivalence, \ie\ there exist a
monoidal functor $J_\mathcal{C}\colon\mathcal{M}^H\to\mathcal{C}$ and
monoidal natural isomorphisms $\eta\colon 1_\mathcal{C}\Rightarrow
I_\mathcal{C}J_\mathcal{C}$ and $\epsilon\colon
J_\mathcal{C}I_\mathcal{C}\Rightarrow 1_\mathcal{D}$.

The inverse of $[I_\mathcal{C},\id_\omega]$ as a functor over
$\Vect_\K$ is given by
$[J_\mathcal{C},\theta]\colon(\mathcal{M}^H,\omega^H)\to(\mathcal{C},\omega)$
where
$\theta:=(\id_{I_\mathcal{C}}\cdot\id_\omega)\circ(\epsilon^{-1}\cdot\id_{\omega^H})\colon\omega^H\Rightarrow
J_\mathcal{C}\cdot\omega$ (Lemma 2.1.3 of~\cite{Sc92}). Obviously,
$\theta$ is monoidal, too. Both $[I_\mathcal{C},\id_\omega]$ and
$[J_\mathcal{C},\theta]$ are thus morphisms of
$\mathfrak{C}^{\otimes,\ast,s}_{\K,rec}$ and mutually inverse.
\end{proof}

\begin{remark}
Tannaka--Kre\v\i n reconstruction is here done for comodules of
coalgebras rather than for modules of algebras for the usual two
reasons.

First, if the coproduct $\Delta\colon C\to C\otimes C$ of the
coalgebra $C$ over $\K$ uses the algebraic tensor product, then the
category of all right-$C$-comodules is already determined by the
finite-dimensional right-$C$-comodules. The forgetful functor
$\omega^C\colon\mathcal{M}^C\to\Vect_\K$ from the category
$\mathcal{M}^C$ of finite-dimensional right-$C$-comodules therefore
maps into $\fdVect_\K$, the full subcategory of $\Vect_\K$ that
contains the rigid objects (those that have left-duals).

Second, when one reconstructs the bialgebra structure of
$\coend(\mathcal{C},\omega)$, one wants the tensor product of the
underlying monoidal category into which $\omega$ maps, to preserve
arbitrary colimits in order to preserve the
$\coend(\mathcal{C},\omega)$, too. This is possible in $\Vect_\K$. In
order to provide an algebra which is reconstructed as a universal
$\mathbf{end}$ with a coalgebra structure, one would need the tensor
product to preserve arbitrary limits. This is not always possible.
\end{remark}

\end{document}